# ASYMPTOTIC BEHAVIOR OF SMOOTH SOLUTIONS
# FOR PARTIALLY DISSIPATIVE HYPERBOLIC SYSTEMS
# WITH A CONVEX ENTROPY

S. BIANCHINI[†], B. HANOUZET[*] AND R. NATALINI[†]

ABSTRACT. We study the asymptotic time behavior of global smooth solutions to general entropy dissipative hyperbolic systems of balance law in $m$ space dimensions, under the Shizuta-Kawashima condition. We show that these solutions approach constant equilibrium state in the $L^p$-norm at a rate $O(t^{-\frac{m}{2}(1-\frac{1}{p})})$, as $t \to \infty$, for $p \in [\min\{m,2\}, \infty]$. Moreover, we can show that we can approximate, with a faster order of convergence, the conservative part of the solution in terms of the linearized hyperbolic operator for $m \geq 2$, and by a parabolic equation, in the spirit of Chapman-Enskog expansion. The main tool is given by a detailed analysis of the Green function for the linearized problem.

1. INTRODUCTION

In the following we shall consider the Cauchy problem for a general hyperbolic symmetrizable $m$-dimensional system of balance laws

$$(1.1) \qquad u_t + \sum_{\alpha=1}^{m} (f_\alpha(u))_{x_\alpha} = g(u),$$

with the initial conditions

$$(1.2) \qquad u(x,0) = u_0(x),$$

where $u = (u_1, u_2) \in \Omega \subseteq \mathbb{R}^{n_1} \times \mathbb{R}^{n_2}$, with $n_1 + n_2 = n$. We also assume that there are $n_1$ conservation laws in the system, namely that we can take

$$(1.3) \qquad g(u) = \begin{pmatrix} 0 \\ q(u) \end{pmatrix}, \text{ with } q(u) \in \mathbb{R}^{n_2}.$$

According to the general theory of hyperbolic systems of balance laws [8], if the flux functions $f_\alpha$ and the source term $g$ are smooth enough, it is well-known that problem (1.1)-(1.2) has a unique local smooth solution, at least for some time interval $[0, T)$ with $T > 0$, if the initial data are also sufficiently smooth. In the general case, and even for very good initial data, smooth solutions may break down in finite time, due to the appearance of singularities, either discontinuities or blow-up.

Despite these general considerations, sometimes dissipative mechanisms due to the source term can prevent the formation of singularities, at least for some restricted classes of initial data, as observed for many models which arise to describe physical phenomena. A typical and well-known example is given by the compressible Euler equations with damping, see [30, 15] for the 1-dimensional case and [34] for an interesting 3-dimensional extension.

Recently, in [13], it was proposed a quite general framework of sufficient conditions which guarantee the global existence in time of smooth solutions. Actually, for the systems which are endowed with a strictly convex entropy function $\mathcal{E} = \mathcal{E}(u)$, a first natural assumption is the *entropy dissipation condition*, see [5, 27, 31, 37], namely for every $u, \overline{u} \in \Omega$, with $g(\overline{u}) = 0$,

$$(\mathcal{E}'(u) - \mathcal{E}'(\overline{u})) \cdot g(u) \leq 0,$$

where $\mathcal{E}'(u)$ is considered as a vector in $\mathbb{R}^n$ and $\cdot$ is the scalar product in the same space. Unfortunately, it is easy to see that this condition is too weak to prevent the formation of singularities, see again [13].

A quite natural supplementary condition can be imposed to entropy dissipative systems, following the classical approach by Shizuta and Kawashima [19, 33], and in the following called condition **(SK)**, which in the present case reads

$$(1.4) \qquad Ker\, Dg(\overline{u}) \cap \{\text{eigenspaces of } \sum_{\alpha=1}^{m} Df_\alpha(\overline{u})\xi_\alpha\} = \{0\},$$

for every $\xi \in \mathbb{R}^m \setminus \{0\}$ and every $\overline{u} \in \Omega$, with $g(\overline{u}) = 0$. It is possible to prove that this condition, which is satisfied in many interesting examples, is also sufficient to establish a general result of global existence for small perturbations of equilibrium constant states, see [13] for a first proof in one space dimension, and [38] for a multidimensional extension.

In this paper we investigate the asymptotic behavior in time of the global solutions, then always assuming the existence of a strictly convex entropy and the **(SK)** condition. First let us describe our approach in the one dimensional case. Our starting point is a careful and refined analysis of the behavior of the Green function for the linearized problem, which is decomposed in three main terms. The first term, the diffusive one, consists of heat kernels moving along the characteristic directions of the local relaxed hyperbolic problem; the singular part consists of exponentially decaying $\delta$-functions along the characteristic directions of the full system. Finally the remainder term decays faster than the first one. Notice that in the one dimensional case a first analysis of the Green function was already contained in [39]. However, there are some differences in our analysis, and in particular we are able to give a more precise description of the behavior of the diffusive part, which is decomposed in four blocks, which decay with different decay rates. This description will be crucial in the analysis of the asymptotic behavior. Let us better explain this point. Actually, we show that solutions have canonical projections on two different components: the conservative part and the dissipative part. The first one, which loosely speaking corresponds to the conservative part of equations in (1.1), decays in time like the heat kernel, since it corresponds to the diffusive part of the Green function. On the other side, the dissipative part is strongly influenced by the dissipation and decays at a rate $t^{-\frac{1}{2}}$ faster of the conservative one. To establish this result, we shall use the Duhamel principle and the Green kernel estimates, very much in the spirit of the Kawashima approach for the hyperbolic-parabolic equations [19]. Unfortunately, with respect to that result, here there is a severe obstruction given by the lack of decay of the source term, when convoluted with the Green kernel. This is not the case for convective or diffusive terms, since they are derivative terms, so having a better decay, of an order $\frac{1}{2}$ for every derivative. However, in the present work we have shown that there is a structural algebraic compatibility between the Green kernel and the conservative structure of the system, by decomposing the kernel according to different linear projectors, which yields the *cancellation* of its highest order and slowly decaying interactions with the source term.

In the multidimensional case, the explicit form of the Green function cannot in general be expressed, and we have to relay directly on the Fourier coordinates. Thus the separation of the Green kernel into various part is done at the level of solution operator $\Gamma(t)$ acting on $L^2(R^m, R^n)$, or $L^1 \cap L^2(R^m, R^n)$. This allow to perform $L^p$ linear decay estimates, for $p \geq 2$.

Let us now shortly review some previous results concerning the asymptotic behavior of solution to dissipative hyperbolic systems. A huge amount of work has been done during many years around the special case of the dissipative (nonlinear) wave equation, see for instance [14, 24, 29] and references therein. At the same time, and starting from the seminal paper by T.P. Liu [20], there were some studies on $2 \times 2$ systems with relaxation, see for instance [15] for the p-system with damping, and [6] for the general case. For more general models, we recall the paper by Y. Zeng [39] about gas dynamics in thermal nonequilibrium and finally the paper by T. Ruggeri and D. Serre [32] about stability of constant equilibrium states for general hyperbolic systems in one space dimension, under zero-mass perturbations. A related result has been recently established by J.F. Coulombel and T. Goudon, who have considered the diffusive relaxation limit of multidimensional isothermal Euler equations [7], see Example 5.12 for a comparison with our approach. Finally let us remark that, under similar assumptions, stability of shock profiles for general relaxation models has been considered in [26].

The paper is organized as follows: Section 2 is devoted to recall some basic results about hyperbolic systems with entropy dissipation and the Shizuta-Kawashima condition. In this section we also introduce the decomposition of the linearized system, which will be called the Conservative-Dissipative form. Section 3 contains a very detailed analysis of the Green kernel in one dimension, while the multidimensional case is presented in Section 4. Finally, Section 5 is devoted to the study of the decay properties of the nonlinear system. Not only we shall prove the decay results for both the conservative and the dissipative part of the solution, but we shall show also that the conservative variable approaches the conservative part of the solution of the corresponding linearized problem, faster that the decay of the heat kernel for $m \geq 2$. Then, we prove that the solution of the parabolic problem, given by the Chapman-Enskog expansion, approximates the conservative part of the solution of the nonlinear hyperbolic system. For $m \geq 2$ the Chapman-Enskog operator is linear, while, in one space dimension, the decay of the nonlinear part has a stronger influence, and so we can only show the faster convergence towards the solution of a parabolic equation with quadratic nonlinearity.

Finally, let us point out again that these results were obtained by assuming all the time the condition **(SK)**. Unfortunately, this condition is not satisfied by many models, as for instance in one space dimension for the Kerr–Debye system, which describes the propagation of electromagnetic waves in nonlinear Kerr medium [12, 16, 13], for perturbations around a null electric field. Another interesting example is given by the equations of gas dynamics in thermal nonequilibrium, which has been investigated in [39], where however global existence of solutions has been established, even if condition **(SK)** is in general violated, thanks to a splitting of the system in two parts, one of them being linearly degenerated. The situation is even worst in more space dimensions, since there are more possibilities to violate condition **(SK)**. This is the case for for every equilibrium state for the 3-dimensional version of Kerr–Debye model, as shown by a simple check. However, a physically relevant class of systems which verify the **(SK)** condition is given by the rotationally invariant systems in Example 4.7 below. Other examples are given by the BGK models proposed in [1], under the Bouchut stability condition [3]. Actually, we expect that, for many physical systems not satisfying the **(SK)** condition, we could consider the influence of other factors, like the existence of linearly degenerate fields, or, in several space dimensions, the well-known faster time decay of the Green function, even for the nondissipative case. Some preliminary results about systems violating the condition **(SK)** will be presented in [25].

## 2. Basic structure of the entropy dissipative hyperbolic systems

### 2.1. Entropy dissipation.

In the following we shall consider a general $m$-dimensional system of balance laws given by equation (1.1), with the source term $g = g(u)$ verifying (1.3).

According to the general theory of hyperbolic systems of conservation laws [8, 31], we shall assume that the system satisfies an entropy principle: there exists a strictly convex function $\mathcal{E} = \mathcal{E}(u)$, the entropy density, and some related entropy-flux functions $\mathcal{F}_\alpha = \mathcal{F}_\alpha(u)$, such that for every smooth solution $u \in \Omega$ to system (1.1), there holds

$$\mathcal{E}_t(u) + \sum_{\alpha=1}^{m} (\mathcal{F}_\alpha(u))_{x_\alpha} = \mathcal{G}(u), \tag{2.1}$$

where $\mathcal{F}_\alpha' = (F_\alpha')^T \mathcal{E}'$ and $\mathcal{G} = \mathcal{E}' \cdot g$. Let us introduce the set $\gamma$ of the equilibrium points:

$$\gamma = \{u \in \Omega; g(u) = 0\}.$$

**Definition 2.1.** *The system (1.1) is* non-degenerate *if, for every $\bar{u} \in \gamma$, it holds*

$$q_{u_2}(\bar{u}) \text{ is non singular}. \tag{2.2}$$

**Definition 2.2.** *The system (1.1) is* entropy dissipative, *if, for every $\bar{u} \in \gamma$ and $u \in \Omega$, we have*

$$(\mathcal{E}'(u) - \mathcal{E}'(\bar{u})) \cdot g(u) \leq 0. \tag{2.3}$$

Following [11, 10, 2], it is now useful to symmetrize our system by introducing a new variable, the *entropy variable*, which is just given by

$$W = \mathcal{E}'(u). \tag{2.4}$$

Actually, since $\mathcal{E}$ is a strictly convex function, we can inverse $\mathcal{E}'$ to recover the original variable $u$ by the inverse map $\Phi \doteq (\mathcal{E}')^{-1}$. Let us set now $A_0(W) \doteq \Phi'(W)$, $C_\alpha(W) \doteq Df_\alpha(\Phi(W))A_0(W)$, and $G(W) = g(\Phi(W)) = \begin{pmatrix} 0 \\ Q(W) \end{pmatrix}$. It is easy to see that the matrix $A_0(W)$ is symmetric positive definite and, for every $\alpha = 1, \ldots, m$, $C_\alpha(W)$ is symmetric. Then, selecting $W$ as the new variable, our system reads

$$(2.5) \qquad A_0(W)W_t + \sum_{\alpha=1}^m C_\alpha(W)W_{x_\alpha} = G(W).$$

Now, as proved in [13], if the system is entropy dissipative and non-degenerate, the set of equilibrium points is locally reduced to a single smooth manifold. More precisely, in the entropy coordinates, we have that, setting $\tilde{\gamma} = \mathcal{E}'(\gamma)$, the entropy dissipation condition reads

$$(2.6) \qquad (W_2 - \bar{W}_2) \cdot Q(W) \leq 0,$$

for every $\bar{W}_2 \in \mathbb{R}^{n_2}$ such that there exists $\bar{W}_1 \in \mathbb{R}^{n_1}$ with $\bar{W} = (\bar{W}_1, \bar{W}_2) \in \tilde{\gamma}$, and every $W \in \mathcal{E}'(\Omega)$. In this case, i.e. if the system is entropy dissipative and non-degenerate, we have that, if $\bar{W} \in \tilde{\gamma}$, then every $W \in \mathcal{E}'(\Omega)$ is also an equilibrium point if and only if $W_2 = \bar{W}_2$, see [13].

Observe now that our definition of dissipative entropy is just invariant for affine perturbations of the form $\tilde{\mathcal{E}}(u) = \mathcal{E}(u) + \alpha + \beta \cdot u$, for $\alpha \in \mathbb{R}$, $\beta \in \mathbb{R}^n$. Therefore, without loss of generality, we can suppose $\bar{u} = 0 \in \gamma$ and consider system (1.1) with $g(0) = 0$, and fix $f_\alpha(0) = 0$. Moreover, we always can assume that the entropy function $\mathcal{E}$ is a quadratic, i.e. such that

$$\mathcal{E}(0) = 0, \ \mathcal{E}'(0) = 0 \in \tilde{\gamma}.$$

Next, following the above considerations, and according to the actual structure of many systems arising in physical models [36, 28, 37, 13, 31], we focus our investigation on a slightly restricted class of entropy dissipative non-degenerate systems, namely the systems such that

$$(2.7) \qquad Q(W) = D(W)W_2, \text{ with } D(0) \text{ negative definite}.$$

In the following we shall refer to these systems just as strictly entropy dissipative systems.

2.2. **The Shizuta-Kawashima condition and the global existence of solutions.** To continue our analysis of smooth solutions for dissipative hyperbolic systems, we need some supplementary coupling conditions to avoid shock formation. A very natural condition was first introduced by Shizuta and Kawashima in [33], for hyperbolic–parabolic systems. Here we first state the condition for the original unknown, i.e. for system (1.1), just assuming that $u = 0$ is an equilibrium point with $g(0) = 0$.

**Definition 2.3.** *The system (1.1) verifies condition* **(SK)**, *if every eigenvector of $\sum_{\alpha=1}^m Df_\alpha(0)\xi_\alpha$ is not in the null space of $Dg(0)$, for every $\xi \in \mathbb{R}^m \setminus \{0\}$.*

Since this condition is invariant under diffeomorphisms which conserve the origin, in the case of strictly entropy dissipative systems, Definition 2.3 is equivalent to

$$(2.8) \qquad \begin{aligned} &\text{for every } \lambda \in \mathbb{R} \text{ and every } X \in \mathbb{R}^{n_1} \setminus \{0\}, \text{ the vector } \begin{pmatrix} X \\ 0 \end{pmatrix} \in \mathbb{R}^n \text{ is not in} \\ &\text{the null space of } \lambda A_0(0) + \sum_{\alpha=1}^m C_\alpha(0)\xi_\alpha, \text{ for every } \xi \in \mathbb{R}^m \setminus \{0\}. \end{aligned}$$

Let us consider now the linearized version of system (1.1), namely, setting $A_\alpha = Df_\alpha(0)$ and $B = Dg(0)$,

$$(2.9) \qquad u_t + \sum_{\alpha=1}^m A_\alpha u_{x_\alpha} = Bu, \qquad u \in \mathbb{R}^n, x \in \mathbb{R}^m, t \in \mathbb{R}^+,$$

with $B$ of the form

$$(2.10) \qquad B = \begin{bmatrix} 0 & 0 \\ D_1 & D_2 \end{bmatrix}, \qquad D_1 \in \mathbb{R}^{n_1 \times n_2}, D_2 \in \mathbb{R}^{n_2 \times n_2},$$

with $n = n_1 + n_2$. Set $A(\xi) = \sum_{\alpha=1}^m A_\alpha \xi_\alpha$. According to the previous discussion, we can assume that (H1) there is a symmetric positive definite matrix $A_0$ such that $A_\alpha A_0$ is symmetric, for every $\alpha = 1, \ldots, m$, and

$$BA_0 = \begin{bmatrix} 0 & 0 \\ 0 & D \end{bmatrix},$$

where $D \in \mathbb{R}^{n_2 \times n_2}$ is negative definite;

(H2) any eigenvector of $A(\xi)$ is not in the null space of $B$, for every $\xi \in \mathbb{R}^m \setminus \{0\}$.

To use the condition **(SK)**, we have to give a reformulation which takes into account the kernel

$$E(i\xi) = B - iA(\xi).$$

This is the content of the following lemma, which is an extension of Theorem 1.1 in [33] to the case of a non symmetric matrix $D$ (the proof is omitted).

**Lemma 2.4.** *Under the assumption (H1), assumption (H2) is equivalent to any of the following:*

i) *there exists $K = K(\xi) \in \mathbb{R}^{n \times n}$ such that, for every $\xi \in \mathbb{R}^m \setminus \{0\}$, $K(\xi)A_0$ is a skew symmetric matrix and*

$$\frac{1}{2}(K(\xi)A(\xi)A_0 + A(\xi)A_0 K^T(\xi)) - \frac{1}{2}(BA_0 + A_0 B^T)$$

*is strictly positive definite;*

ii) *if $\lambda(z)$ is an eigenvalue of $E(z)$, then $\Re(\lambda(i\xi)) < 0$ for every $\xi \in \mathbb{R}^m \setminus \{0\}$;*

iii) *there exists $c > 0$ such that*

(2.11) $$\Re(\lambda(i\xi)) \leq -c \frac{|\xi|^2}{1 + |\xi|^2},$$

*for every $\xi \in \mathbb{R}^m \setminus \{0\}$.*

About the existence of a solution, we recall the following result [13, 38].

**Theorem 2.5.** *Assume that system (1.1) is strictly entropy dissipative and condition **(SK)** is satisfied. Then there exists $\delta > 0$ such that, if $\|u_0\|_s \leq \delta$, with $s \geq [m/2] + 2$, there is a unique global solution $u$ of (1.1)–(1.2), which verifies*

$$u \in C^0([0, \infty); H^s(\mathbb{R}^m)) \cap C^1([0, \infty); H^{s-1}(\mathbb{R}^m)),$$

*and such that, in terms of the entropy variable $W = (W_1, W_2)$,*

(2.12) $$\sup_{0 \leq t < +\infty} \|W(t)\|_s^2 + \int_0^{+\infty} \left(\|\nabla W_1(\tau)\|_{s-1}^2 + \|W_2(\tau)\|_s^2\right) d\tau \leq C(\delta) \|W_0\|_s^2,$$

*where $C(\delta)$ is a positive constant.*

2.3. **The Conservative-Dissipative form in the linear case.** We now consider a linear system with constant coefficients:

(2.13) $$w_t + \sum_{\alpha=1}^m \tilde{A}_\alpha w_{x_\alpha} = \tilde{B} w,$$

where $w = (w_1, w_2) \in \mathbb{R}^{n_1} \times \mathbb{R}^{n_2}$. We assume also that the differential part is symmetric:

(2.14) $$\text{for all } \alpha = 1, \ldots, m, \qquad \tilde{A}_\alpha^T = \tilde{A}_\alpha.$$

**Definition 2.6.** *Under assumption (2.14), the partially dissipative system (2.13) is in Conservative-Dissipative form (C-D form) if there exists a negative definite matrix $\tilde{D} \in \mathbb{R}^{n_2}$, such that*

(2.15) $$\tilde{B} = \begin{pmatrix} 0 & 0 \\ 0 & \tilde{D} \end{pmatrix}.$$

In the following $w_1 \doteq w_c$ is called the conservative variable, while $w_2 \doteq w_d$ is the dissipative one.

We notice that, thanks to assumption (H1), system (2.9) is already in the C-D form if $A_0 = I$. We shall prove in the following that there exists a linear change of variable such that (2.9) takes the C-D form in the general case of $A_0$ symmetric and positive definite.

Take $u$ a solution of (2.9). First, we use the classical transformation

$$v = A_0^{-1/2} u,$$

which yields

$$(2.16) \qquad v_t + \sum_{\alpha=1}^{m} \bar{A}_\alpha v_{x_\alpha} = \bar{B}v,$$

where $\bar{A}_\alpha = A_0^{-1/2} A_\alpha A_0^{1/2}$, $\bar{B} = A_0^{-1/2} B A_0^{1/2}$. Notice that system (2.16) has a symmetric differential part, but the matrix $\bar{B}$ does not satisfy (2.15). However, by the assumptions on the matrix $B$, for $\bar{B}$ there exists a null space of dimension $n_1$, while the other eigenvalues are strictly negative. We shall construct the C-D variables using the projection $Q_0$ on the null space and the complementing projection $Q_- = I - Q_0$. We compute $Q_0$ by using the explicit formula, see [18]:

$$Q_0 = -\frac{1}{2\pi i} \oint_{|\xi| \ll 1} (\bar{B} - \xi I)^{-1} d\xi.$$

We have:

$$(\bar{B} - \xi I)^{-1} = A_0^{1/2} (BA_0 - \xi A_0)^{-1} A_0^{1/2} = A_0^{1/2} \begin{bmatrix} -\xi A_{0,11} & -\xi A_{0,21} \\ -\xi A_{0,21} & D - \xi A_{0,22} \end{bmatrix}^{-1} A_0^{1/2}$$

$$= A_0^{1/2} \left( \begin{bmatrix} -\xi A_{0,11} & -\xi A_{0,12} \\ 0 & D \end{bmatrix} (I + O(\xi)) \right)^{-1} A_0^{1/2}$$

$$= A_0^{1/2} (I + O(\xi)) \begin{bmatrix} -(A_{0,11})^{-1}/\xi & -(A_{0,11})^{-1} A_{0,12} D^{-1} \\ 0 & D^{-1} \end{bmatrix} A_0^{1/2}$$

$$= \frac{1}{\xi} A_0^{1/2} \begin{bmatrix} -(A_{0,11})^{-1} & 0 \\ 0 & 0 \end{bmatrix} A_0^{1/2} + O(1).$$

We thus obtain that

$$Q_0 = A_0^{1/2} \begin{bmatrix} (A_{0,11})^{-1} & 0 \\ 0 & 0 \end{bmatrix} A_0^{1/2}.$$

Note that due to the assumptions on $A_0$, this projector is symmetric. In particular we can choose left and right projectors $L_0 \in \mathbb{R}^{n \times n_1}$, $R_0 \in \mathbb{R}^{n_1 \times n}$, so that

$$(2.17) \qquad Q_0 = R_0 L_0, \qquad L_0 R_0 = I \in \mathbb{R}^{n_1 \times n_1}, \qquad L_0 = R_0^T.$$

Note that by the last condition also $R_0, L_0$ are unique: in fact they are given by

$$(2.18) \qquad R_0 = A_0^{1/2} \begin{bmatrix} (A_{0,11})^{-1/2} \\ 0 \end{bmatrix}, \qquad L_0 = \begin{bmatrix} (A_{0,11})^{-1/2} & 0 \end{bmatrix} A_0^{1/2}.$$

We define the complementary projection $Q_-$ to be

$$(2.19) \qquad Q_- \doteq I - Q_0 = R_- L_-, \qquad L_- R_- = I \in \mathbb{R}^{n_2 \times n_2}, \qquad L_- = R_-^T.$$

The last condition follows because also $Q_-$ is symmetric, and the matrices $R_- \in \mathbb{R}^{n_2 \times n}$, $L_- \in \mathbb{R}^{n \times n_2}$ are the unique left and right projectors which satisfy (2.19): one can check that these projectors are given by

$$(2.20) \qquad R_- = A_0^{-1/2} \begin{bmatrix} 0 \\ ((A_0^{-1})_{22})^{-1/2} \end{bmatrix}, \qquad L_- = \begin{bmatrix} 0 & ((A_0^{-1})_{22})^{-1/2} \end{bmatrix} A_0^{-1/2}.$$

Set now $w_1 = L_0 v$, $w_2 = L_- v$. We have

$$v = (Q_0 + Q_-)v = R_0 L_0 v + R_- L_- v$$

$$= R_0 w_1 + R_- w_2$$

and by (2.16):

$$(L_0 v)_t + \sum_{\alpha=1}^{m} L_0 \bar{A}_\alpha (R_0 w_1 + R_- w_2)_{x_\alpha} = L_0 \bar{B}(R_0 w_1 + R_- w_2),$$

$$(L_- v)_t + \sum_{\alpha=1}^{m} L_- \bar{A}_\alpha (R_0 w_1 + R_- w_2)_{x_\alpha} = L_- \bar{B}(R_0 w_1 + R_- w_2).$$

Now notice that
$$L_0 \bar{B} = 0, \quad \bar{B} R_0 = 0.$$

Therefore $w = (w_1, w_2)$ are Conservative-Dissipative variables and system (2.16) is equivalent to the C-D form system (2.13), where $\tilde{A}_\alpha$ are the symmetric matrices

(2.21) $$\tilde{A}_\alpha = \begin{pmatrix} L_0 \bar{A}_\alpha R_0 & L_0 \bar{A}_\alpha R_- \\ L_- \bar{A}_\alpha R_0 & L_- \bar{A}_\alpha R_- \end{pmatrix}$$

and

(2.22) $$\tilde{B} = \begin{pmatrix} 0 & 0 \\ 0 & \tilde{D} \end{pmatrix},$$

with

(2.23) $$\tilde{D} = L_- \bar{B} R_- = ((A_0^{-1})_{22})^{1/2} D ((A_0^{-1})_{22})^{1/2}$$

is negative definite.

**Proposition 2.7.** *If $u$ is a solution to system* (2.9), *then, under assumption (H1)*,

(2.24) $$w = Mu = \begin{bmatrix} (A_{0,11})^{-1/2} & 0 \\ ((A_0^{-1})_{22})^{-1/2}(A_0^{-1})_{21} & ((A_0^{-1})_{22})^{1/2} \end{bmatrix} u$$

*is a solution to the C-D form system* (2.13) *with* (2.21), (2.22) *and* (2.23).

*Remark* 2.8. If system (1.1) is non degenerate and entropy-dissipative, we can apply Proposition 2.7 to the linearized system (2.9). Therefore, in the following, we are always going to assume that the unknown $u$ is chosen in such a way that (2.9) is in conservative-dissipative form. In this case, we say that also system (1.1) is in conservative-dissipative form and we shall set $u = (u_c, u_d) \in \mathbb{R}^{n_1} \times \mathbb{R}^{n_2}$.

*Remark* 2.9. More generally, we can look to the set of linear transformations $w = Mu$, such that, starting from system (2.9), under assumption (H1), the new system is in C-D form. To obtain the symmetry of the differential part, we have to take $M$ such that $(M^T M)^{-1}$ is a symmetrizer of the system. Hence, we can choose $M$ such that

(2.25) $$M^T M = A_0^{-1}.$$

Now, to verify condition (2.15), we obtain the relations

(2.26) $$\begin{cases} M_{11}^T M_{11} = (A_{0,11})^{-1}, \\ M_{12} = 0, \\ M_{21} = (M_{22}^T)^{-1}(A_0^{-1})_{21}, \\ M_{22}^T M_{22} = (A_0^{-1})_{22}. \end{cases}$$

In particular, a special choice is to take $M_{11}$ and $M_{22}$ symmetric and we obtain (2.24) with $\tilde{D}$ given by (2.23).

*Example* 2.10. **The $p$-system with relaxation.** Let us consider system

(2.27) $$\begin{cases} \partial_t u + \partial_x v = 0, \\ \partial_t v + \partial_x \sigma(u) = h(u) - v, \end{cases}$$

with $\sigma'(u) > 0$. Its linear counterpart is given by

(2.28) $$\begin{cases} \partial_t u + \partial_x v = 0, \\ \partial_t v + \lambda^2 \partial_x u = au - v, \end{cases}$$

where $\lambda = \sqrt{\sigma'(0)}$ and $a = h'(0)$. Therefore

$$A = \begin{pmatrix} 0 & 1 \\ \lambda^2 & 0 \end{pmatrix}, \quad B = \begin{pmatrix} 0 & 0 \\ a & -1 \end{pmatrix}.$$

We can use the symmetrizer $A_0$ given by

$$A_0 = \begin{pmatrix} 1 & a \\ a & \lambda^2 \end{pmatrix},$$

which is positive definite if it holds the subcharacteristic condition $\lambda > |a|$. It is easy to verify that assumption (H1) is verified, since

$$AA_0 = \begin{pmatrix} a & \lambda^2 \\ \lambda^2 & a\lambda^2 \end{pmatrix}, \quad BA_0 = \begin{pmatrix} 0 & 0 \\ 0 & a^2 - \lambda^2 \end{pmatrix}.$$

To recover the C-D form, we first compute the inverse matrix $A_0^{-1}$, which is given by

$$A_0^{-1} = \frac{1}{\lambda^2 - a^2} \begin{pmatrix} \lambda^2 & -a \\ -a & 1 \end{pmatrix}.$$

This yields

$$M = \begin{pmatrix} 1 & 0 \\ -a(\lambda^2 - a^2)^{-\frac{1}{2}} & (\lambda^2 - a^2)^{-\frac{1}{2}} \end{pmatrix}$$

and so we obtain the matrices of the C-D form

$$\tilde{A} = \begin{pmatrix} a & (\lambda^2 - a^2)^{\frac{1}{2}} \\ (\lambda^2 - a^2)^{\frac{1}{2}} & -a \end{pmatrix}, \quad \tilde{B} = \begin{pmatrix} 0 & 0 \\ 0 & -1 \end{pmatrix}.$$

Setting

$$\begin{pmatrix} u_c \\ u_d \end{pmatrix} = M \begin{pmatrix} u \\ v \end{pmatrix}$$

and reporting in (2.27), we obtain its conservative-dissipative form

(2.29) $$\partial_t \begin{pmatrix} u_c \\ u_d \end{pmatrix} + \partial_x \begin{pmatrix} au_c + (\lambda^2 - a^2)^{\frac{1}{2}} u_d \\ (\lambda^2 - a^2)^{-\frac{1}{2}}(\sigma(u_c) - a^2 u_c) - au_d \end{pmatrix} = \begin{pmatrix} 0 \\ (\lambda^2 - a^2)^{-\frac{1}{2}}(h(u_c) - au_c) - u_d \end{pmatrix}.$$

3. THE GREEN KERNEL FOR LINEAR DISSIPATIVE SYSTEMS IN ONE SPACE DIMENSION

Aim of this section is to compute the Green kernel $\Gamma(t)$ for a linear dissipative hyperbolic system. The fact that we are in dimension one will help us in inverting the Fourier transforms, hence giving explicit form to the principal parts of $\Gamma(t)$.

We can consider directly a system in C-D form, according to the results of Subsection 2.3. So we write our system as

(3.1) $$w_t + Aw_x = Bw,$$

where $w = (w_c, w_d) \in \mathbb{R}^{n_1} \times \mathbb{R}^{n_2}$. We assume that the differential part is symmetric, and there exists a negative definite matrix $D \in \mathbb{R}^{n_2 \times n_2}$, such that

(3.2) $$B = \begin{pmatrix} 0 & 0 \\ 0 & D \end{pmatrix},$$

so we have (H1). We assume also that (3.1) verifies condition **(SK)**, then we have also (H2). We notice that, by contrast with [39], we are not assuming that the matrix $B$ is symmetric. This is necessary to deal with some specific examples, as for instance the Jin-Xin relaxation system, see [17, 13].

We want to study the Green kernel $\Gamma(t, x)$ of (3.1), which satisfies

(3.3)
$$\begin{cases} \Gamma_t + A\Gamma_x = B\Gamma \\ \Gamma(0, x) = \delta(x)I \end{cases}$$

Taking the Fourier transform

$$\hat{\Gamma}(t, \xi) = \int_{\mathbb{R}} \Gamma(t, x)e^{-i\xi x}dx$$

of (3.3) we obtain

(3.4)
$$\begin{cases} d\hat{\Gamma}/dt = (B - i\xi A)\hat{\Gamma} \\ \hat{\Gamma}(0, \xi) = I \end{cases}$$

To study the large time behavior of the Green kernel $\Gamma$, we use the approach already proposed in [21], [39].

3.1. **Perturbation analysis.** Consider the entire function

(3.5)
$$E(z) = B - zA.$$

It is clear that the solution to (3.4) is given by

(3.6)
$$\hat{\Gamma}(t, \xi) = e^{E(i\xi)t} = \sum_{n=0}^{+\infty} \frac{t^n}{n!}(B - i\xi A)^n,$$

so that

(3.7)
$$\Gamma(t, x) = \text{p.v.} \frac{1}{2\pi} \int_{\mathbb{R}} e^{E(i\xi)t}e^{i\xi x}d\xi = \lim_{N\to+\infty} \frac{1}{2\pi} \int_{-N}^{N} e^{E(i\xi)t}e^{i\xi x}d\xi$$

and $\hat{\Gamma}(t, z) = e^{E(z)t}$ is an entire function of $z$. The next analysis follows using some ideas in [18].

The function $E(z)$ given by (3.5), as a matrix valued function, has a constant number $s$ of distinct eigenvalues $\lambda(z)$ iff $z$ is not one of the exceptional points, which are of finite number in the plane. In fact these points are the solutions to

(3.8)
$$\det(B - zA - \lambda I),$$

which is a polynomial equation with holomorphic coefficients. It follows that its roots $\lambda(z)$ are branches of one or more analytic functions with algebraic singularities of at most order $n$. As a consequence the number of eigenvalues is constant, with the exception of a finite number of points, called *exceptional points*, in each compact set of the complex plane. Since we can write

$$E(z) = z(-A + B/z),$$

then the same occurs in a neighborhood of $z = \infty$, so that in our case the number of exceptional points is bounded in the whole complex plane. Even if $z$ is not an exceptional point, differently from [39], the matrix $E(z)$ is in general not diagonalizable, due to the fact that $B$ is negative definite but not symmetric: we say that $E(z)$ is *permanently degenerate*.

In any region where there are no exceptional point, the functions $\lambda_j(z)$, $j = 1, \ldots, s$, are holomorphic, with constant multiplicities $m_j$, $j = 1, \ldots, s$. In general these $\lambda_j$ are branches of one or more algebraic functions, denoted again as $\lambda_j$, $j = 1, \ldots, s$. The exceptional points can be either regular points for these algebraic functions, or a branch point for some $\lambda_j(z)$. In the first case the eigenprojectors remain bounded, while in a branch point the projectors have a pole.

In general, the function $E(z)$, if $z$ is not exceptional, is represented as

(3.9)
$$E(z) = \sum_j \lambda_j(z)P_j(z) + \sum_j D_j(z),$$

where $\lambda_j$ are the eigenvalues of $E(z)$, $P_j(z)$ the corresponding eigenprojections, given by the formula

(3.10)
$$P_j(z) = -\frac{1}{2\pi i} \oint_{|\xi - \lambda_j(z)| \ll 1} (E(z) - \xi I)^{-1}d\xi,$$

and $D_j$ are the nilpotent matrices, due to the fact that in general $E$ is not diagonalizable, defined by

$$D_j(z) = (E(z) - \lambda_j(z)I)P_j(z). \tag{3.11}$$

Note that by construction the eigenvalues of $D_j$ are 0, so that

$$D_j^{m_j}(z) = 0, \tag{3.12}$$

where $m_j$ is the multiplicity of $\lambda_j$.

We now study which consequences have the assumptions (H1), (H2) on $E(z)$ and $\hat{\Gamma}(t,z)$ near the point $z = 0$ and $z = \infty$. Both points are in general exceptional points: for $z \to 0$, $n_1$ eigenvalues different from 0 converges to 0. When $|z| \to \infty$, the matrix $A$ is diagonalizable, but it can have common eigenvalues: then the perturbation $B/z$ will in general remove part of this degeneracy.

We are going to show that, near $z = 0$, semisimple eigenvalue of $B$, the matrix $E(z)$ has a decomposition

$$E(z) = \sum_{jk} \left( \Lambda_{jk}(z) P_{jk}(z) + D_{jk}(z) \right) + E_1(z), \tag{3.13}$$

where the $\Lambda_{jk}$ are diagonal $n \times n$ matrices composed by the $n_1$ eigenvalues, which converge to 0, the $P_{jk}$ are spectral projectors, the $D_{jk}$ are nilpotent operators commuting with $P_{jk}$. From assumption (H1), we can control the behavior of the matrix $E_1(z)$.

In a similar way, near $z = \infty$, since the eigenvalues of $A$ are semisimple, $E(z)$ has a canonical decomposition as

$$E(z) = \sum_{jk} \left( \Upsilon_{jk}(z) \mathcal{P}_{jk}(z) + \mathcal{D}_{jk}(z) \right). \tag{3.14}$$

The entries in the matrix $\Lambda_{jk}$ have an expansion in the form

$$\lambda(z) = -z\lambda_j^1 - z^2 c_{jk} + \mathcal{O}(z^3), \tag{3.15}$$

where the $\lambda_j^1$ are the eigenvalues of the symmetric block $A_{11}$.

On the other hand, the entries in the matrix $\Upsilon_{jk}$ have an expansion in the form

$$v(z) = -z\lambda_j + b_{jk} + \mathcal{O}(1/z), \tag{3.16}$$

where the $\lambda_j$ are the eigenvalues of $A$. As a consequence of the assumption (H2), which is equivalent to the condition (**SK**), the coefficients $c_{jk}$ and $b_{jk}$ have strictly negative real part.

3.1.1. *Case $z = 0$.* The total projector $P$ corresponding to all the eigenvalues near 0 is

$$P(z) = -\frac{1}{2\pi i} \oint_{|\xi| \ll 1} (E(z) - \xi I)^{-1} d\xi. \tag{3.17}$$

The point $z = 0$ is in general an exceptional point, and the projections corresponding to the eigenvalues with negative real part (i.e. not in any of the $jk$ families defined in (3.15)) can have poles in $z = 0$. Nevertheless, the projection

$$P_-(z) = I - P(z) \tag{3.18}$$

corresponding to the whole family of eigenvalues with strictly negative real part is holomorphic near $z = 0$, see [18] or the analysis below.

To simplify computations, we introduce the projectors

$$L_0 = R_0^T = \begin{bmatrix} I_{n_1} & 0 \end{bmatrix}, \quad L_- = R_-^T = \begin{bmatrix} 0 & I_{n_2} \end{bmatrix}. \tag{3.19}$$

and

$$Q_0 = R_0 L_0 = \begin{bmatrix} I_{n_1} & 0 \\ 0 & 0 \end{bmatrix}, \quad Q_- = I - Q_0 = \begin{bmatrix} 0 & 0 \\ 0 & I_{n_2} \end{bmatrix}. \tag{3.20}$$

For $\xi$ close to 0, we have

$$(B - \xi I)^{-1} = -\frac{1}{\xi} Q_0 + R_-(D - \xi I)^{-1} L_- = \begin{bmatrix} \xi^{-1} I & 0 \\ 0 & (D - \xi I)^{-1} \end{bmatrix}. \tag{3.21}$$

We recall also the expansion of the resolvent

$$R(\xi, z) \doteq (B - zA - \xi I)^{-1} = (B - \xi I)^{-1}\left(I - zA(B - \xi I)^{-1}\right)^{-1}$$
$$= (B - \xi I)^{-1} + z(B - \xi I)^{-1}A(B - \xi I)^{-1} + z^2(B - \xi I)^{-1}(A(B - \xi I)^{-1})^2 + O(z^3)$$
(3.22)
$$= R_0(\xi) + \sum_{n \geq 1} z^n R_n(\xi).$$

The total projector $P$ becomes here

(3.23)
$$P(z) = P_0 + \sum_{n \geq 1} z^n P_n,$$

where $P_n$ is given by the integral

(3.24)
$$P_n = -\frac{1}{2\pi i} \oint_{|\xi| \ll 1} R_n(\xi) d\xi.$$

By using (3.21) we have the zero order coefficient,

(3.25)
$$P_0 = -\frac{1}{2\pi i} \oint_{|\xi| \ll 1} (B - \xi I)^{-1} d\xi = Q_0 = \begin{bmatrix} I & 0 \\ 0 & 0 \end{bmatrix},$$

while the coefficient for $z$ is given by

$$P_1 = -\frac{1}{2\pi i} \oint_{|\xi| \ll 1} (B - \xi I)^{-1} A (B - \xi I)^{-1} d\xi = Q_0 A R_- D^{-1} L_- + R_- D^{-1} L_- A Q_0$$
(3.26)
$$= \begin{bmatrix} 0 & A_{12} D^{-1} \\ D^{-1} A_{21} & 0 \end{bmatrix}.$$

For completeness we will also compute the coefficient for $z^2$. Integrating as before $R_2(\xi)$, we have

$$P_2 = Q_0 (A R_- D^{-1} L_-)^2 + R_- D^{-1} L_- A Q_0 A R_- D^{-1} L_- + (R_- D^{-1} L_- A)^2 Q_0$$
$$- (Q_0 A)^2 R_- D^{-2} L_- - R_- D^{-2} L_- (A Q_0)^2 - Q_0 A R_- D^{-2} L_- A Q_0$$
(3.27)
$$= \begin{bmatrix} -A_{12} D^{-2} A_{21} & A_{12} D^{-1} A_{22} D^{-1} - A_{11} A_{12} D^{-2} \\ D^{-1} A_{22} D^{-1} A_{21} - D^{-2} A_{21} A_{11} & D^{-1} A_{21} A_{12} D^{-1} \end{bmatrix}.$$

As we will see later, the coefficient we are interested is the 22 coefficient: in fact we see that we can write

(3.28)
$$P(z) = \begin{bmatrix} I + O(z^2) & z A_{12} D^{-1} + O(z^2) \\ z D^{-1} A_{21} + O(z^2) & z^2 D^{-1} A_{21} A_{12} D^{-1} + O(z^3) \end{bmatrix}.$$

We introduce the right and left eigenprojectors of $P(z)$, $R(z) \in \mathbb{R}^{n \times n_1}$, $L(z) \in \mathbb{R}^{n_1 \times n}$, which verify

$$P(z) = R(z) L(z), \quad L(z) R(z) = I.$$

We can find the power series of $L(z)$ and $R(z)$ by means of the relations

$$L(z) P(z) = L(z), \qquad P(z) R(z) = R(z).$$

We have in fact for $L(z) = L_0 + z L_1 + z^2 L_1 + O(z^3)$, with $L_0$ given by (3.19),

$$\left(L_0 + z L_1 + z^2 L_2 + O(z^3)\right)\left(P_0 + z P_1 + z^2 P_2 + O(z^3)\right) = L_0 P_0 + z\left(L_0 P_1 + L_1 P_0\right)$$
$$+ z^2\left(L_0 P_2 + L_1 P_1 + L_2 P_0\right) + O(z^3)$$
$$= L_0 + z L_0 A R_- D^{-1} L_- + z^2 L_0 (A R_- D^{-1} L_-)^2 - z^2 L_0 A R_- D^{-2} L_- A Q_0$$
$$- z^2 L_0 A Q_0 A R_- D^{-2} L_- + z^2 L_1 R_- D^{-1} L_- A Q_0 + z^2 L_2 P_0 + O(z^3),$$

so that we see that

(3.29)
$$L_1 = L_0 A R_- D^{-1} L_- = \begin{bmatrix} 0 & A_{12} D^{-1} \end{bmatrix}, \quad R_1 = R_- D^{-1} L_- A R_0 = \begin{bmatrix} 0 \\ D^{-1} A_{21} \end{bmatrix},$$

$$L_2 = L_0(AR_-D^{-1}L_-)^2 - L_0AQ_0AR_-D^{-2}L_- = \begin{bmatrix} -A_{12}D^{-2}A_{21}/2 & A_{12}D^{-1}A_{22}D^{-1} - A_{11}A_{12}D^{-2} \end{bmatrix},$$

(3.30)
$$R_2 = (R_-D^{-1}L_-A)^2R_0 - R_-D^{-2}L_-AQ_0AR_0 = \begin{bmatrix} -A_{12}D^{-2}A_{21}/2 \\ D^{-1}A_{22}D^{-1}A_{21} - D^{-2}A_{21}A_{11} \end{bmatrix},$$

where we used a similar computation for $R_1$, $R_2$. Note that since $D^{-1}$ is not symmetric, these projectors do not satisfy $L(z) = R(z)$.

Next, we can decompose $E(z)$ according to the right and left operators:

(3.31) $$E(z) = R(z)F(z)L(z) + R_-(z)F_-(z)L_-(z),$$

where $F(z) \doteq L(z)E(z)R(z) \in \mathbb{R}^{n_1 \times n_1}$ and $F_-(z) \doteq L_-(z)E(z)R_-(z) \in \mathbb{R}^{n_2 \times n_2}$.

We have
$$F(z) = \left(L_0 + zL_1 + O(z^2)\right)(B - zA)\left(R_0 + zR_1 + O(z^2)\right)$$
$$= -zL_0AR_0 - z^2L_0AR_-D^{-1}L_-AR_0 + O(z^3)$$
(3.32)
$$= -zA_{11} - z^2A_{12}D^{-1}A_{21} + O(z^3).$$

The matrix $A_{11}$ is symmetric, from assumption (H1), so that we can write for some eigenvalues $\lambda_j^1$, with multiplicity $m'_j$, $j = 1, \ldots, m'$, and left and right eigenprojections $l_j \in \mathbb{R}^{m'_j \times n_1}$, $r_j \in \mathbb{R}^{n_1 \times m'_j}$, with $l_j = r_j^T$,

(3.33) $$A_{11} = \sum_{j=1}^{m'} \lambda_j^1 r_j l_j.$$

**Lemma 3.1.** *Under the assumption (H2), the matrix $A_{12}D^{-1}A_{21}$ is negative defined.*

*Proof.* Let $r$ be an eigenvector of $A_{11}$ for the eigenvalue $\lambda^1$. Then, since $D^{-1}$ is strictly negative and $A_{12} = A_{21}^T$, we have that
$$d = r^T A_{12}D^{-1}A_{21}r < 0$$
if $A_{21}r \neq 0$. On the other hand, we have that
$$B\begin{bmatrix} r \\ 0 \end{bmatrix} = 0$$
and
$$A\begin{bmatrix} r \\ 0 \end{bmatrix} = A\begin{bmatrix} A_{11}r \\ A_{21}r \end{bmatrix} = \begin{bmatrix} \lambda^1 r \\ A_{21}r \end{bmatrix}.$$

Therefore $\begin{bmatrix} r \\ 0 \end{bmatrix}$ is an eigenvalue of $A$ if $A_{21}r = 0$, and assumption (H2) implies that $A_{21}r \neq 0$. □

We can again reduce (3.32) by considering the right and left projections $r_j(z)$, $l_j(z)$, with $r_j(0) = r_j$, $l_j(0) = l_j$, for each family of eigenvalues $\lambda_j = -z\lambda_j^1 + O(z^2)$. We can now expand the projectors as
$$p_j(z) = r_j l_j + z p_j^1 + O(z^2), \quad r_j(z) = r_j + z r_j^1 + O(z^2), \quad l_j(z) = l_j + z l_j^1 + O(z^2),$$
where $r_j^1 = p_j^1 r_j$, $l_j^1 = p_j^1 l_j$. Therefore
$$F_j(z) \doteq l_j(z)F(z)r_j(z) = -z\lambda_j^1 I_{m_j} - z^2\left(l_j A_{12}D^{-1}A_{21}r_j\right) - z^2\left(l_j^1 A_{11}r_j + l_j A_{11}r_j^1\right) + O(z^3).$$

Now, using formula (**II**-2.14) in [18], it is possible to prove that the third term vanishes. So, we obtain

(3.34) $$F_j(z) = -z\lambda_j^1 I - z^2(A_{21}r_j)^T D^{-1}(A_{21}r_j) + O(z^3).$$

As before, $(A_{21}r_j)^T D^{-1}(A_{21}r_j)$ has eigenvalues with strictly negative real part. Let $c_{jk}$ be the eigenvalues, with multiplicity $m'_j$, of the reduced matrix $(A_{21}r_j)^T D^{-1}(A_{21}r_j)$, and let $p_{jk} \in \mathbb{R}^{m'_j \times m'_j}$ be the corresponding eigenprojections, with $d_{jk} \in \mathbb{R}^{m'_j \times m'_j}$ the nilpotent matrices. We obtain finally that

(3.35) $$F_j(z) = \sum_k F_{jk}(z) = \sum_k (-z\lambda_j^1 - z^2 c_{jk})p_{jk} - z^2 d_{jk} + O(z^3).$$

The eigenvalues tending to 0 belong to $jk$ families, whose $z$ expansion can be expressed by

$$\lambda_{jk}(z) = -z\lambda_j^1 - z^2 c_{jk} + O(z^3). \tag{3.36}$$

Note that since all $c_{jk}$ are different, then the total projection for each family is holomorphic near $z = 0$, and similarly the nilpotent part. They can be expressed as

$$P_{jk}(z) = R_0 r_j p_{jk}(R_0 r_j)^T + O(z), \qquad D_{jk}(z) = z^2 R_0 r_j d_{jk}(R_0 r_j)^T + O(z^3). \tag{3.37}$$

Therefore, we obtain the projection of $E(z)$ on the null eigenvalue as

$$R(z)F(z)L(z) = \sum_{jk} \left(-z\lambda_j^1 I - z^2 c_{jk} I + O(z^3)\right) P_{jk}(z) + D_{jk}(z). \tag{3.38}$$

*Remark* 3.2. As we have noticed before, in general inside the $jk$ family there are different eigenvalues whose projections have a pole in $z = 0$. We just say that the total projection $P_{jk}(z)$ of the whole $jk$ family does not have poles in 0: as we showed, this follows because $F(z)$ can be decomposed as the sum of $F_{jk}(z)$, acting on different subspaces.

We study now the term $F_-(z) = L_-(z)E(z)R_-(z)$. As before, we expand the left and right projections of $P_-(z)$ as $L_-(z) = L_- + zL_-^1 + O(z^2)$ and $R_-(z) = R_- + zR_-^1 + O(z^2)$. We obtain

$$L_-^1 = \begin{bmatrix} -D^{-1}A_{21} & 0 \end{bmatrix}, \qquad R_-^1 = \begin{bmatrix} -A_{12}D^{-1} \\ 0 \end{bmatrix}. \tag{3.39}$$

This yields

$$F_-(z) = D - zA_{22} + O(z^2). \tag{3.40}$$

We sum up the previous results in the following statement.

**Proposition 3.3.** *We have the following decomposition near $z = 0$*

$$E(z) = \sum_{jk} \left(\Lambda_{jk}(z)P_{jk}(z) + D_{jk}(z)\right) + E_1(z), \tag{3.41}$$

*where the $\Lambda_{jk}$ are diagonal $n \times n$ matrices composed by the $n_1$ eigenvalues $\lambda_{jk}$ given by (3.36), the coefficients $c_{jk}$ having strictly negative real part, thanks to assumption (H2). The spectral projectors $P_{jk}$ and the nilpotent operators $D_{jk}$ are given by (3.37) and verify*

$$P_{jk}(z)P_{j'k'}(z) = \delta_{jj'}\delta_{kk'}P_{jk}(z), \quad D_{jk}(z)P_{jk}(z) = P_{jk}(z)D_{jk}(z) = D_{jk}(z),$$

$$\Lambda_{jk}(z)P_{jk}(z) = P_{jk}(z)\Lambda_{jk}(z) = \Lambda_{jk}(z), \quad P_{jk}(z)E_1(z) = E_1(z)P_{jk}(z) = 0.$$

*The term $E_1(z)$ is given by $R_-(z)F_-(z)L_-(z)$, where $F_-(0) = D$, which, by assumption (H1), has eigenvalues with strictly negative real part.*

3.1.2. *Case $z = \infty$.* We do now the same analysis when $|z| \to \infty$. We have

$$E(z) = B - zA = z\left(-A + \frac{1}{z}B\right) = z\tilde{E}(1/z),$$

where

$$\tilde{E}(\zeta) = -A + \zeta B. \tag{3.42}$$

Since $A$ is symmetric, we can write

$$A = \sum_j \lambda_j R_j R_j^T, \tag{3.43}$$

where $\lambda_j$ are the eigenvalues with multiplicity $m_j$, $R_j \in \mathbb{R}^{n \times m_j}$ are the right eigenprojections, normalized by $R_j^T R_j = I \in \mathbb{R}^{m_j \times m_j}$. As before, by considering the total projection for the family of eigenvalues converging to $\lambda_j$ as $\zeta \to 0$, we obtain the reduced equation for each $\lambda_j$,

$$\tilde{F}_j(\zeta) = -\lambda_j I + \zeta R_j^T B R_j + O(\zeta^2). \tag{3.44}$$

If one decompose now $R_j^T B R_j$ as

$$R_j^T B R_j = \sum_k b_{jk} \tilde{p}_{jk} + \tilde{d}_{jk}, \qquad \tilde{p}_{jk}, \tilde{d}_{jk} \in \mathbb{R}^{m_{jk} \times m_{jk}},$$

we obtain that we can reduce further the $\tilde{F}_j$ by

$$\tilde{F}_{jk}(\zeta) = (-\lambda_j + \zeta b_{jk})\tilde{p}_{jk} + \zeta \tilde{d}_{jk} + O(\zeta^2). \tag{3.45}$$

As before, one obtains the $jk$ families of eigenvalues for $|z| \to \infty$ have the $z$ series

$$\lambda_{jk}(z) = -z\lambda_j + b_{jk} + O(1/z), \tag{3.46}$$

and the projectors and nilpotent parts

$$\mathcal{P}_{jk} = R_j \tilde{p}_{jk} R_j^T + O(1/z), \qquad \mathcal{D}_{jk} = R_j \tilde{d}_{jk} R_j^T + O(1/z). \tag{3.47}$$

**Lemma 3.4.** *Under the assumption (H2), the eigenvalues $b_{jk}$ of $R_j^T B R_j$ have strictly negative real part.*

The proof follows by arguing as in Lemma 3.1.

**Proposition 3.5.** *We have the following decomposition near $z = \infty$*

$$E(z) = \sum_{jk} \big(\Upsilon_{jk}(z)\mathcal{P}_{jk}(z) + \mathcal{D}_{jk}(z)\big), \tag{3.48}$$

*where $\Upsilon_{jk}$ is the diagonal matrix whose entries are the eigenvalues of the $jk$ family (3.46), the coefficients $b_{jk}$ having strictly negative real part, thanks to assumption (H2). The spectral projectors $\mathcal{P}_{jk}$ and the nilpotent operators $\mathcal{D}_{jk}$ are given by (3.47) and verify*

$$\mathcal{P}_{jk}(z)\mathcal{P}_{j'k'}(z) = \delta_{jj'}\delta_{kk'}\mathcal{P}_{jk}(z), \quad \mathcal{D}_{jk}(z)\mathcal{P}_{jk}(z) = \mathcal{P}_{jk}(z)\mathcal{D}_{jk}(z) = \mathcal{D}_{jk}(z),$$

$$\Upsilon_{jk}(z)\mathcal{P}_{jk}(z) = \mathcal{P}_{jk}(z)\Upsilon_{jk}(z) = \Upsilon_{jk}(z).$$

**3.2. Green function estimates.** In the general case, assuming that the matrix $A$ is symmetric, we have that $\Gamma(t,x) = 0$ if $x > \bar{\lambda} t$ or $x < \underline{\lambda} t$, where

$$\bar{\lambda} := \max_j \lambda_j, \quad \underline{\lambda} := \min_j \lambda_j, \tag{3.49}$$

namely, the support of $\Gamma$ is contained in the wave cone of $A$. Therefore we have

$$\Gamma(t,x) = \Gamma(t,x)\chi\{\underline{\lambda} t \leq x \leq \bar{\lambda} t\}, \tag{3.50}$$

where $\chi$ is the characteristic function. In conclusion, in the following we shall assume all the time

$$\left|\frac{x}{t}\right| \leq C.$$

Now, we are ready to estimate the global behavior for large $t$ of the Green kernel $\Gamma(t,x)$ using the local expansions contained in Propositions 3.3 and 3.5. We associate a diffusive operator with Green function $K(t,x)$ to the expansion (3.41), and a dissipative transport operator with Green function $\mathcal{K}(t,x)$ to (3.48), and we estimate the remainder term

$$R(t,x) = \Gamma(t,x) - K(t,x) - \mathcal{K}(t,x).$$

*3.2.1. Estimates near $z = 0$.* In the following we shall consider the Green kernel as composed of 4 parts, acting on $w_c$, $w_d$:

$$\Gamma(t,x) = \begin{bmatrix} L_0 \Gamma(t,x) R_0 & L_0 \Gamma(t,x) R_- \\ L_- \Gamma(t,x) R_0 & L_- \Gamma(t,x) R_- \end{bmatrix} = \begin{bmatrix} \Gamma_{00}(t,x) & \Gamma_{0-}(t,x) \\ \Gamma_{-0}(t,x) & \Gamma_{--}(t,x) \end{bmatrix}. \tag{3.51}$$

Using the expansion of $L(z)$, $R(z)$ given by (3.29), (3.30) it follows that

$$P_{jk}(z) = \bar{P}_{jk}(z) + R_{jk}(z)$$

$$= \begin{bmatrix} r_j p_{jk} r_j^T & z r_j p_{jk} r_j^T A_{12} D^{-1} \\ z D^{-1} A_{21} r_j p_{jk} r_j^T & z^2 D^{-1} A_{21} r_j p_{jk} r_j^T A_{12} D^{-1} \end{bmatrix} + \begin{bmatrix} O(z^2) & O(z^2) \\ O(z^2) & O(z^3) \end{bmatrix}. \tag{3.52}$$

We can associate to each term of $F_{jk}(i\xi)$ given by (3.35) the the parabolic equation

$$(3.53) \qquad w_t + \lambda_j^1 w_x = -(c_{jk}I + d_{jk})w_{xx}, \qquad w \in \mathbb{R}^{m'_j},$$

where $c_{jk}$ is in general complex valued, but its real part is strictly negative:

$$(3.54) \qquad c_{jk} = -\mu_{jk} - i\nu_{jk}, \qquad \mu_{jk} > 0.$$

Its kernel can be computed explicitly. If $\gamma_{jk} = \sqrt{-c_{jk}} = \sqrt{\mu_{jk} + i\nu_{jk}}$ is the square root with positive real part, so that $\arg \gamma_{jk} \in (-\pi/4, \pi/4)$, then

$$(3.55) \qquad g_{jk}(t,x) \doteq \frac{1}{2\gamma_{jk}\sqrt{\pi t}} \exp\left\{-\frac{(x-\lambda_j^1 t)^2}{4(\mu_{jk}+i\nu_{jk})t}\right\} \left[\sum_\iota M_{jk,\iota} \frac{(x-\lambda_j^1 t)^{2\iota}}{((\mu_{jk}+i\nu_{jk})t)^\iota}\right].$$

The matrix valued coefficients $M_{jk,\iota}$ are due to the fact that we have a nilpotent part $D_{jk}$: the maximal value of $\iota$ is $m_{jk} - 1$, where $m_{jk}$ is the multiplicity of $c_{jk}$. Note that we have in any case that for some $c > 0$

$$(3.56) \qquad |g_{jk}(t,x)| \le \frac{O(1)}{\sqrt{t}} e^{-(x-\lambda_j^1 t)^2/(ct)} \quad \forall k,\ (t,x) \in \mathbb{R}^+ \times \mathbb{R}.$$

Similarly, one can see that the inverse Fourier transform of

$$\bar{F}(t,z) \doteq \sum_{jk} e^{-z\lambda_j^1 t - z^2 c_{jk} t} e^{-z^2 D_{jk} t} \bar{P}_{jk}(z)$$

is given by the function

$$(3.57) \qquad K(t,x) \doteq \sum_{jk} \begin{bmatrix} r_j(g_{jk}(t,x)p_{jk})r_j^T & -r_j\left(\dfrac{dg_{jk}}{dx}p_{jk}\right)r_j^T A_{12} D^{-1} \\ -D^{-1} A_{21} r_j\left(\dfrac{dg_{jk}}{dx}p_{jk}\right)p_{jk}r_j^T & D^{-1} A_{21} r_j\left(\dfrac{d^2 g_{jk}}{dx^2}p_{jk}\right)r_j^T A_{12} D^{-1} \end{bmatrix}.$$

The function $K(t,x)$, as we will see later, collects the principal parts of each component (3.51) of the Green kernel $\Gamma(t,x)$.

By the Proposition 3.3, we can compute $e^{E(z)t}$ near $z = 0$:

$$(3.58) \qquad e^{E(z)t} = \sum_{jk} e^{-z\lambda_j^1 t - z^2 c_{jk} t} e^{-z^2 t D_{jk} + O(z^3 t)} P_{jk}(z) + R_-(z) e^{F_-(z)t} L_-(z).$$

We associate to kernel of the parabolic equation (3.53), the function

$$\hat{g}_{jk}(z) = -(z\lambda_j^1 - z^2 c_{jk})I - z^2 d_{jk}.$$

In the same way, to the Green function $K(t,x)$ we associate the function

$$(3.59) \qquad \hat{K}(z) \doteq \sum_{jk} \hat{g}_{jk}(z) \bar{P}_{jk}(z) = \sum_{jk} \begin{bmatrix} r_j \hat{g}_{jk}(z) p_{jk} r_j^T & z r_j \hat{g}_{jk}(z) p_{jk} r_j^T A_{12} D^{-1} \\ z D^{-1} A_{21} r_j \hat{g}_{jk}(z) p_{jk} r_j^T & z^2 D^{-1} A_{21} r_j \hat{g}_{jk}(z) p_{jk} r_j^T A_{12} D^{-1} \end{bmatrix}.$$

Consider the following integral:

$$R_1(t,x) = \frac{1}{2\pi} \int_{-\epsilon}^{\epsilon} \left(e^{E(i\xi)t} - e^{\hat{K}(i\xi)t}\right) e^{i\xi x} d\xi$$

$$= \sum_{jk} \frac{1}{2\pi} \int_{-\epsilon}^{\epsilon} e^{i\xi(x-\lambda_j^1 t) + \xi^2 c_{jk} t} \left(e^{\xi^2 t D_{jk} + O(\xi^3 t)} P_{jk}(i\xi) - e^{\xi^2 t D_{jk}} \bar{P}_{jk}(i\xi)\right) d\xi$$

$$(3.60) \qquad + \frac{1}{2\pi} \int_{-\epsilon}^{\epsilon} R_-(i\xi) e^{F_-(i\xi)t} L_-(i\xi) e^{i\xi x} d\xi.$$

The constant $\epsilon$ is sufficiently small. The meaning of the above integral is that for low frequencies $\xi$, the main parts of $\Gamma(t,x)$ is given by the parabolic diffusion process described by (3.53). Moreover, we are taking into account the principal parts of each component of $\Gamma(t,x)$, as in (3.51).

Using (3.40), and since $D$ has strictly negative eigenvalues, it is clear that for some positive constant $C$

$$(3.61) \qquad \left|\frac{1}{2\pi} \int_{-\epsilon}^{\epsilon} R_-(i\xi) e^{F_-(i\xi)t} L_-(i\xi) e^{i\xi x} d\xi\right| \le C e^{-t/C}.$$

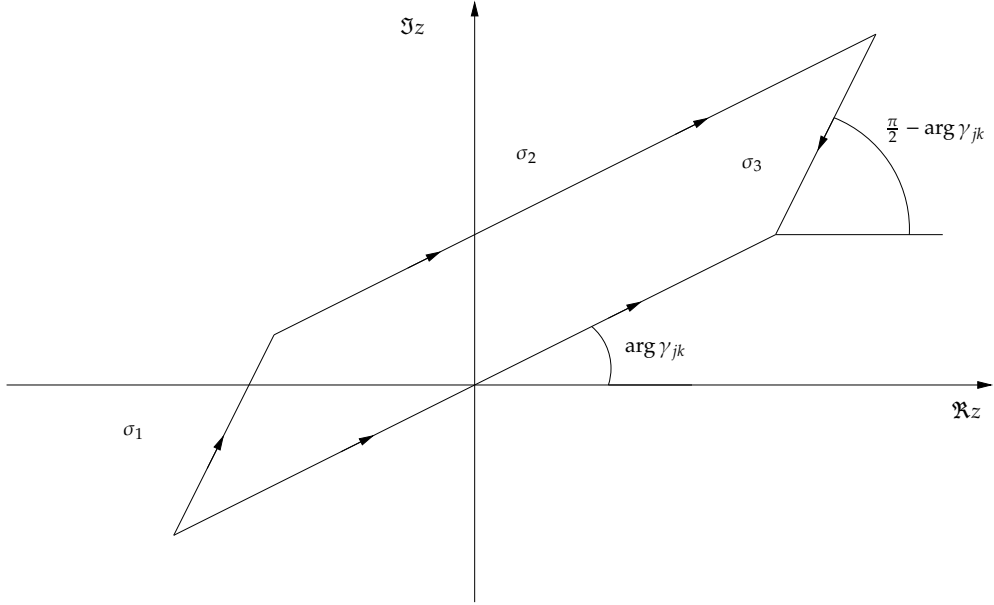

FIGURE 1. The path $\sigma$ in the complex plane in the case $x - \lambda_j^1 t > 0$.

We will consider separately each of the integrals

$$(3.62) \quad R_{jk}(t,x) = \frac{1}{2\pi} \int_{-\epsilon}^{\epsilon} e^{i\xi(x-\lambda_j^1 t) + \xi^2 c_{jk} t} \left( e^{\xi^2 t D_{jk} + O(\xi^3 t)} P_{jk}(i\xi) - e^{\xi^2 t D_{jk}} \bar{P}_{jk}(i\xi) \right) d\xi.$$

Since the integrand is holomorphic (because we are considering the whole eigenspace $P_{jk}$, see Remark 3.2), we can change the path of integration in such a way that, when we take the real part of the exponent inside the integral, we will obtain a strictly negative exponent (it is clear that such a thing do not happens for the path considered in (3.67)). Let $c_{jk} = -\mu_{jk} - i\nu_{jk}$, and denote as before $\gamma_{jk} = \sqrt{\mu_{jk} + i\nu_{jk}}$ its square root with positive real part. Note that for all $jk$ we have

$$(3.63) \quad \arg \gamma_{jk} \in (-\pi/4 + \zeta, \pi/4 - \zeta), \qquad \zeta > 0.$$

By the change of coordinates $\xi = e^{-i \arg \gamma_{jk}} z$, we can write the exponent as

$$(3.64) \quad i\xi(x - \lambda_j^1 t) - \xi^2 \gamma_{jk}^2 t = ie^{-i \arg \gamma_{jk}} z(x - \lambda_j^1 t) - |\gamma_{jk}|^2 z^2 t,$$

integrated along the path $z = e^{i \arg \gamma_{jk}} \xi$, $\xi \in [-\epsilon, \epsilon]$. Since all integrands (3.67) are holomorphic, we can deform the path as in fig. 1: denoting with $y$ the constant

$$(3.65) \quad y = \min \left\{ |x - \lambda_j^1 t|/(2|\gamma_{jk}|^2 t), \epsilon/2 \right\},$$

consider the path

$$\sigma = \left\{ -\epsilon e^{i \arg \gamma_{jk}} + i \operatorname{sgn}(x - \lambda_j^1 t) \eta e^{-i \arg \gamma_{jk}}; \eta \in [0, y] \right\}$$

$$\bigcup \left\{ -\eta e^{i \arg \gamma_{jk}} + i \operatorname{sgn}(x - \lambda_j^1 t) y e^{-i \arg \gamma_{jk}}; \eta \in [-\epsilon, \epsilon] \right\}$$

$$\bigcup \left\{ \epsilon e^{i \arg \gamma_{jk}} - i \operatorname{sgn}(x - \lambda_j^1 t) \eta e^{-i \arg \gamma_{jk}}; \eta \in [-y, 0] \right\}$$

$$(3.66) \quad = \sigma_1 \cup \sigma_2 \cup \sigma_3.$$

We now estimate:

$$(3.67) \quad R_{jk}(t,x) = \frac{1}{2\pi} \oint_{\sigma} e^{ie^{-i \arg \gamma_{jk}} z(x-\lambda_j^1 t) - |\gamma_{jk}|^2 z^2 t} \left( e^{\xi^2 t D_{jk} + O(\xi^3 t)} P_{jk}(i\xi) - e^{\xi^2 t D_{jk}} \bar{P}_{jk}(i\xi) \right) \bigg|_{\xi = e^{-i \arg \gamma_{jk}} z} e^{-i \arg \gamma_{jk}} dz.$$

We begin with the following lemma.

ASYMPTOTIC BEHAVIOR FOR PARTIALLY DISSIPATIVE HYPERBOLIC SYSTEMS17**Lemma 3.6.** *If $D$ is a complex nilpotent matrix, then for all $\delta > 0$ and $\alpha \in \mathbb{C}$, there exists $C = C(\delta)$ such that*

$$\left|e^{\alpha D + A} - e^{\alpha D}\right| \leq C e^{|\alpha|\delta + C|A|}|A| \quad (3.68)$$

*for every matrix $A$.*

*Proof.* Let $0 < \omega << 1$. Since $D$ is nilpotent, there exists an invertible change of base $R = R(\omega)$, with

$$|R| = \mathcal{O}(1), \ |R^{-1}| = \mathcal{O}(\omega^{-n+1}),$$

such that the matrix

$$Y := RDR^{-1} \quad (3.69)$$

verifies

$$|Y| = \omega.$$

In fact, there exists $S \in \mathbb{R}^{n \times n}$, such that

$$D' = SDS^{-1} = \begin{pmatrix} 0 & 1 & 0 & \cdots & \cdots \\ \vdots & 0 & \ddots & \ddots & \\ \vdots & & \ddots & 1 & \ddots \\ & & & \ddots & 0 & \ddots \\ & & & & \ddots & \ddots & 0 \\ & & & & & & 0 \\ \cdots & \cdots & & & \cdots & & 0 \end{pmatrix}.$$

We introduce the diagonal matrix

$$T(\omega) = diag(\omega^{n-1}, \omega^{n-2}, \ldots, \omega, 1)$$

and its inverse

$$T(\omega)^{-1} = diag(\omega^{-n+1}, \omega^{-n+2}, \ldots, \omega^{-1}, 1),$$

and set $R(\omega) = T(\omega)S$, which yields (3.69) with $Y = \omega D'$. We can now compute

$$\begin{aligned}\left|e^{-\alpha D} e^{\alpha D + A} - I\right| &= \left|\left(\sum_{n=0}^{+\infty} \frac{(-\alpha D)^n}{n!}\right)\left(\sum_{m=0}^{+\infty} \frac{(\alpha D + A)^m}{m!}\right) - I\right| \\ &= \left|R^{-1}\left(\sum_{n=1}^{+\infty} \sum_{m=0}^{n-1} \frac{1}{m!(n-m)!}(-\alpha Y)^m (\alpha Y + RAR^{-1})^{n-m}\right)R\right| \\ &\leq C \sum_{n=1}^{+\infty} \sum_{m=0}^{n-1} \frac{1}{m!(n-m)!} \sum_{i=1}^{n-m} \frac{(n-m)!}{i!(n-m-i)!} |\alpha Y|^{n-i} |RAR^{-1}|^i \\ &\leq C \sum_{n=1}^{+\infty} \sum_{i=1}^{n} \frac{1}{i!(n-i)!} |\alpha Y|^{n-i} |RAR^{-1}|^i \sum_{m=0}^{n-i} \frac{(n-i)!}{m!(n-i-m)!} \\ &= C \sum_{n=1}^{+\infty} 2^n \sum_{i=1}^{n} \frac{1}{i!(n-i)!} (|\alpha||Y|)^{n-i} |RAR^{-1}|^i \\ &= C \sum_{n=0}^{+\infty} \frac{2^n}{n!}\left((|\alpha||Y| + |RAR^{-1}|)^n - (|\alpha||Y|)^n\right) \\ &= C\left(e^{2(|\alpha||Y| + |RAR^{-1}|)} - e^{2|\alpha||Y|}\right) \leq C e^{2|\alpha||Y| + C|A|}|A|.\end{aligned}$$

Therefore, we obtain

$$\left|e^{\alpha D + A} - e^{\alpha D}\right| \leq C e^{3|\alpha||Y| + C|A|}|A| \quad (3.70)$$

and the conclusion follows. $\square$

We introduce now

$$\Delta_{jk}(z) = e^{\xi^2 t D_{jk} + \mathcal{O}(\xi^3 t)} - e^{\xi^2 t D_{jk}}.$$

Recall that we set $\xi = e^{-i \arg \gamma_{jk}} z$. Using (3.68), we obtain the estimate

$$\left|\Delta_{jk}(z)\right| \leq C e^{|z|^2 t \delta + \mathcal{O}(|z|^3 t)} |z|^3 t \leq C e^{2|z|^2 t \delta} |z|^3 t, \quad (3.71)$$

with $\epsilon \approx \delta$ and $\delta \ll 1$. Note that we only assume here that $|z|$ is small, but it may happens that $z^3 t$ is large. The importance of the above lemma is in the fact that we do not require any assumption on the norm of $A$.

We have

$$\left|e^{\xi^2 tD_{jk}+O(\xi^3 t)}P_{jk}(i\xi) - e^{\xi^2 tD_{jk}}\bar{P}_{jk}(i\xi)\right| \leq \left\|\begin{bmatrix} r_j \Delta_{jk}(z)p_{jk}r_j^T & zr_j\Delta_{jk}(z)p_{jk}r_j^T A_{12}D^{-1} \\ zD^{-1}A_{21}r_j\Delta_{jk}(z)p_{jk}p_{jk}r_j^T & z^2 D^{-1}A_{21}r_j\Delta_{jk}(z)p_{jk}r_j^T A_{12}D^{-1} \end{bmatrix}\right\|$$

$$+ \left|e^{\xi^2 tD_{jk}+O(\xi^3 t)}\left(P_{jk}(i\xi) - \bar{P}_{jk}(i\xi)\right)\right|$$

Using (3.71), (3.68) and (3.69), we obtain

$$(3.72) \quad \left|e^{\xi^2 tD_{jk}+O(\xi^3 t)}P_{jk}(i\xi) - e^{\xi^2 tD_{jk}}\bar{P}_{jk}(i\xi)\right| = (|z^3|t + |z|)e^{2\delta|z|^2 t}\begin{bmatrix} O(1) & O(1)|z| \\ O(1)|z| & O(1)|z|^2 \end{bmatrix}.$$

What we are going to obtain is the following result: the principal terms of $\Gamma(t, x)$ are the heat kernels $g_{jk}(t, x)$ or their derivatives, and the error terms for each principal part are of higher order. Using (3.67), we will thus integrate along the path $\sigma$ the function

$$(3.73) \quad |z|^p(|z|^2 t + 1)\exp\left\{\Re\left(ie^{-i\arg\gamma_{jk}}z(x - \lambda_j^1 t) - z^2|\gamma_{jk}|^2 t\right) + 2\delta|z|^2 t\right\},$$

for $p = 1, 2, 3$.

On $\sigma_1$, observing that $|z| = O(\epsilon)$, one has
(3.74)
$$\left|\frac{1}{2\pi}\oint_{\sigma_1}(\ldots)dz\right| \leq C\epsilon^p(\epsilon^2 t + 1)\int_0^y \exp\left\{-\cos(2\arg\gamma_{jk})\left(\eta|x - \lambda_j^1 t| + |\gamma_{jk}|^2\epsilon^2 t - |\gamma_{jk}^2|\eta^2 t\right) + 2\delta(\epsilon^2 t)\right\}d\eta$$

$$\leq C\epsilon^{p+1}(\epsilon^2 t + 1)\exp\left\{-\tfrac{1}{2}\mu_{jk}\epsilon^2 t\right\} \leq C\epsilon^{p+1}\exp\left\{-\tfrac{\epsilon^2}{C}t\right\} \leq Ce^{-t/C},$$

for some large constant $C$ and if $\delta$ is sufficiently small. The same estimate can be obtained on $\sigma_3$.

On $\sigma_2$ we have
(3.75)
$$\left|\frac{1}{2\pi}\oint_{\sigma_2}(\ldots)dz\right|$$

$$\leq C\int_{-\epsilon}^{\epsilon}\exp\left\{-\cos(2\arg\gamma_{jk})\left(y|x - \lambda_j^1 t| + |\gamma_{jk}|^2\eta^2 t - |\gamma_{jk}|^2 y^2 t\right) + 2\delta(|z|^2 t)\right\}\left(t(|y| + |\eta|)^2 + 1\right)(|y| + |\eta|)^p d\eta.$$

Recall now that $y = \min\{|x - \lambda_j^1 t|/(2|\gamma_{jk}|^2 t), \epsilon/2\}$. If $\epsilon < |x - \lambda_j^1 t|/|\gamma|^2 t$, then we can evaluate (3.75) as before,

$$(3.76) \quad C\epsilon^{p+1}(\epsilon^2 t + 1)\exp\left\{-\frac{1}{8}\mu_{jk}\epsilon^2 t\right\} \leq Ce^{-t/C}.$$

For $\epsilon \geq |x - \lambda_j^1 t|/|\gamma|^2 t$ we have

(3.77)
$$\left| \frac{1}{2\pi} \oint_{\sigma_2} (\dots) dz \right|$$

$$\leq C \int_{-\epsilon}^{\epsilon} \exp\left\{ -\cos(2 \arg \gamma_{jk}) \left( \frac{(x - \lambda_j^1 t)^2}{4|\gamma_{jk}|^2 t} + |\gamma_{jk}|^2 \eta^2 t \right) + 4\delta \eta^2 t + \frac{\delta(x - \lambda_j^1 t)^2}{|\gamma_{jk}|^4 t} \right\} \left( t(|y| + |\eta|)^2 + 1 \right)(|y| + |\eta|)^p d\eta$$

$$\leq C \exp\left\{ -\cos(2 \arg \gamma_{jk}) \frac{(x - \lambda_j^1 t)^2}{8|\gamma|^2 t} \right\} \int_{-\epsilon}^{\epsilon} e^{-\frac{1}{2} \mu_{jk} \eta^2 t} \left( |\eta|^2 t + y^2 t + 1 \right)(|\eta| + y)^p d\eta$$

$$\leq C \exp\left\{ -\cos(2 \arg \gamma_{jk}) \frac{(x - \lambda_j^1 t)^2}{8|\gamma|^2 t} \right\} \frac{1}{t^{(1+p)/2}} \left( \sum_{\iota=0}^{p+2} \left( \frac{|x - \lambda_j^1 t|}{\sqrt{t}} \right)^{\iota} \right)$$

$$\leq \frac{C}{t^{(1+p)/2}} e^{-(x - \lambda_j^1 t)^2/ct}.$$

Observe here that in (3.62), $R_{jk}$ is bounded for $0 \leq t \leq 1$. Then, using also (3.50), we can write that the rest part near $\xi = 0$ is of the order of

(3.78) $$R_1(t, x) = \sum_j e^{-(x - \lambda_j^1 t)^2/ct} \begin{bmatrix} O(1)(1 + t)^{-1} & O(1)(1 + t)^{-3/2} \\ O(1)(1 + t)^{-3/2} & O(1)(1 + t)^{-2} \end{bmatrix}.$$

3.2.2. *Estimates near $z = \infty$.* In this case, we can associate to each term of $z\tilde{F}_{jk}(1/z)$ in (3.45), the Fourier transform of the Green kernel of the transport equation

(3.79) $$w_t + \lambda_j w_x = (b_{jk} I + \tilde{d}_{jk}) w, \qquad w \in \mathbb{R}^{m_{jk}}.$$

We can write it explicitly by

(3.80) $$\tilde{g}_{jk}(t, x) = \delta(x - \lambda_j t) e^{b_{jk} t} \sum_{\iota} \frac{t^{\iota}}{\iota!} (\tilde{d}_{jk})^{\iota},$$

Note that by our assumption $\mathfrak{R}(b_{jk}) \leq -c < 0$, so that we have the estimate

(3.81) $$|\tilde{g}_{jk}(t, x)| \leq C \delta(x - \lambda_j t) e^{-ct} \quad \forall k, (t, x) \in \mathbb{R}^+ \times \mathbb{R}.$$

We associate to the kernels $\tilde{g}_{jk}$, the hyperbolic Green function

(3.82) $$\mathcal{K}(t, x) = \sum_{jk} R_j \tilde{g}_{jk}(t, x) \tilde{p}_{jk} R_j^T.$$

Using (3.47) and Proposition 3.5, we obtain

(3.83) $$\mathcal{K}(t, x) = \sum_{jk} \delta(x - \lambda_j t) e^{b_{jk} t} e^{t \mathcal{D}_{jk}(\infty)} \mathcal{P}_{jk}(\infty).$$

Observe that $\mathcal{K}(t, x)$ is the Fourier transform in $x$ of

(3.84) $$\hat{\mathcal{K}}(t, \xi) = \sum_{jk} e^{-i \lambda_j t \xi + b_{jk} t} e^{t \mathcal{D}_{jk}(\infty)} \mathcal{P}_{jk}(\infty).$$

From the Proposition 3.5, near $z = \infty$ we have the expansion for $E(z)$ as

(3.85) $$E(z) = \sum_{jk} \left( \left( -z \lambda_j + b_{jk} + b_{jk}^1 \frac{1}{z} + O(1) \frac{1}{z^2} \right) I + \mathcal{D}_{jk}(z) \right) \mathcal{P}_{jk}(z),$$

so that we can compute its exponential as

(3.86) $$e^{E(z) t} = \sum_{jk} e^{-z \lambda_j t + b_{jk} t + (b_{jk}^1 \frac{1}{z} + O(1) \frac{1}{z^2}) t} e^{t \mathcal{D}_{jk}(z)} \mathcal{P}_{jk}(z).$$

For $N$ sufficiently large, we define

$$(3.87) \qquad R_2(t,x) = \frac{1}{2\pi}\int_{|\xi|\geq N}\left(e^{E(i\xi)t}-\hat{\mathcal{K}}(t,\xi)\right)e^{i\xi x}d\xi.$$

We have

$$|R_2(t,x)|\leq\sum_{jk}\frac{1}{2\pi}\left|\int_{|\xi|\geq N}e^{i\xi(x-\lambda_j t)+b_{jk}t}\left(e^{t\mathcal{D}_{jk}(i\xi)+t(b_{jk}^1\frac{1}{i\xi}+O(1)\frac{1}{\xi^2})I}\mathcal{P}_{jk}(i\xi)-e^{t\mathcal{D}_{jk}(\infty)}\mathcal{P}_{jk}(\infty)\right)d\xi\right|.$$

We have

$$\mathcal{D}_{jk}(i\xi)+b_{jk}^1\frac{1}{i\xi}I+O(1)\frac{1}{\xi^2}=\mathcal{D}_{jk}(\infty)+\frac{1}{\xi}\mathcal{D}_{jk}^1+O(1)\left(\frac{1}{\xi^2}\right),\quad \mathcal{P}_{jk}(i\xi)=\mathcal{P}_{jk}(\infty)+\frac{1}{\xi}\mathcal{P}_{jk}^1+O(1)\left(\frac{1}{\xi^2}\right).$$

We denote by $de^D$ the derivative of the application

$$D \to e^D.$$

So, it holds

$$|e^{D+A}-e^D-de^D A|\leq C|A|^2\sup_{s\in[0,1]}|d^2 e^{D+sA}|\leq Ce^{|D|+|A|}|A|^2,$$

which in the case

$$D\mapsto tD,\quad |D|\leq\delta,\quad A\mapsto\left(\frac{t}{\xi}\right)A,\quad |A|/|\xi|\leq\delta,$$

becomes

$$\left|e^{tD+(\frac{t}{\xi})A}-e^{tD}-de^{tD}\left(\frac{t}{\xi}\right)A\right|\leq Ce^{t(|D|+|(\frac{1}{\xi})A|)}\left(\frac{t^2}{\xi^2}\right)|A|^2\leq C\left(\frac{t^2}{\xi^2}\right)e^{C\delta t}.$$

From the Lemma 3.6 and (3.69) we can suppose that $|\mathcal{D}_{jk}(\infty)|\leq\delta$.

Therefore, for $N$ large enough, we obtain

$$e^{t\mathcal{D}_{jk}(i\xi)+t(b_{jk}^1\frac{1}{i\xi}+O(1)\frac{1}{\xi^2})I}\mathcal{P}_{jk}(i\xi) = \left[e^{t\mathcal{D}_{jk}(\infty)}\mathcal{P}_{jk}(\infty)-d\left(e^{t\mathcal{D}_{jk}(\infty)}\right)\left(\frac{t}{\xi}\right)\left(\mathcal{D}_{jk}^1+O(1)\left(\frac{1}{\xi}\right)\right)+O(1)\left(\frac{t^2}{\xi^2}\right)e^{C\delta t}\right]$$

$$\times\left[\mathcal{P}_{jk}(\infty)+\frac{1}{\xi}\mathcal{P}_{jk}^1+O(1)\left(\frac{1}{\xi^2}\right)\right].$$

We conclude that

$$(3.88)\quad e^{t\mathcal{D}_{jk}(i\xi)+t(b_{jk}^1\frac{1}{i\xi}+O(1)\frac{1}{\xi^2})I}\mathcal{P}_{jk}(i\xi)-e^{t\mathcal{D}_{jk}(\infty)}\mathcal{P}_{jk}(\infty)=\frac{1}{\xi}\left[td\left(e^{t\mathcal{D}_{jk}(\infty)}\right)\mathcal{D}_{jk}^1\mathcal{P}_{jk}(\infty)+e^{t\mathcal{D}_{jk}(\infty)}\mathcal{P}_{jk}^1\right]+\frac{1}{\xi^2}O(1)e^{C\delta t}.$$

Finally, for $\delta$ small enough, we obtain

$$(3.89)\qquad |R_2(t,x)|\leq\sum_{jk}Ce^{-\alpha t}\left|\int_{|\xi|\geq N}\frac{e^{i\xi(x-\lambda_j t)}}{\xi}d\xi\right|+\sum_{jk}Ce^{-\alpha t}\int_{|\xi|\geq N}\frac{1}{\xi^2}d\xi\leq Ce^{-\alpha t}.$$

3.2.3. *Estimates in between.* To complete the study of the Fourier transform of $\Gamma$, we have to study which terms are left: these are the parabolic kernel $K$ for $|\xi|\geq\epsilon$ and $t\geq 1$, the transport kernel $\mathcal{K}$ for $|\xi|\leq N$, and the kernel $E(z)$ for $\epsilon\leq|\xi|\leq N$. We thus have to consider here 3 cases.

First, one has immediately that if $K$ is the parabolic linearized Green kernel. Set

$$(3.90)\qquad R_3(t,x)=\frac{1}{2\pi}\int_{|\xi|\geq\epsilon}\hat{K}(t,\xi)e^{i\xi x}d\xi,$$

and we obtain

$$(3.91)\qquad |R_3(t,x)|\leq\sum_{jk}C\left|\int_{|\xi|\geq\epsilon}e^{i\xi(x-\lambda_j^1 t)}e^{-c_{jk}\xi^2 t}d\xi\right|\leq C\frac{e^{-t/C}}{\sqrt{t}},\qquad t\geq 1.$$

Similarly, we introduce

$$(3.92)\qquad R_4(t,x)=\frac{1}{2\pi}\int_{|\xi|\leq N}\hat{\mathcal{K}}(t,\xi)e^{i\xi x}d\xi,$$

$$(3.93) \quad |R_4(t,x)| \leq Ce^{-\alpha t}\sum_j \left|\int_{-N}^{N} e^{i\xi(x-\lambda_j t)}d\xi\right| \leq Ce^{-\alpha t}\min_j\{N, 1/|x-\lambda_j^1 t|\}.$$

Finally, we set

$$(3.94) \quad R_5(t,x) = \frac{1}{2\pi}\int_{\epsilon \leq |\xi| \leq N} e^{E(i\xi)t}e^{i\xi x}d\xi,$$

and, thanks to Lemma 2.4, we can use the estimate (2.11), which follows from the Shizuta-Kawashima condition:

$$(3.95) \quad |R_5(t,x)| \leq C\int_{\epsilon \leq |\xi| \leq N} e^{-\alpha\xi^2 t/(1+\xi^2)}d\xi \leq Ce^{-t/C},$$

for some large constant $C$.

3.2.4. *Global estimates.* Notice that, thanks to (3.50), we have to study only the case $|x/t| \leq C$. By means of (3.78), (3.89), (3.91), (3.93), (3.95), we have

$$(3.96) \quad \begin{aligned} R(t,x) &= \Gamma(t,x) - K(t,x) - \mathcal{K}(t,x) \\ &= \sum_j \frac{e^{-(x-\lambda_j^1 t)^2/Ct}}{1+t}\begin{bmatrix} O(1) & O(1)(1+t)^{-1/2} \\ O(1)(1+t)^{-1/2} & O(1)(1+t)^{-1} \end{bmatrix} + Ce^{-t/C}, \end{aligned}$$

for some large constant $C$. Moreover, using again $|x/t| \leq C$, we have that in (3.96) the first term on the RHS dominates the second one, so that we can write

$$(3.97) \quad R(t,x) = \sum_j \frac{e^{-(x-\lambda_j^1 t)^2/Ct}}{1+t}\begin{bmatrix} O(1) & O(1)(1+t)^{-1/2} \\ O(1)(1+t)^{-1/2} & O(1)(1+t)^{-1} \end{bmatrix}.$$

We have then proved the main result of this section.

**Theorem 3.7.** *Let $\Gamma(t,x)$ be the Green function of system (3.1), under the assumptions (H1) and (H2). Let $K(t,x)$ be the Green function of the diffusive operator given by (3.57) and $\mathcal{K}(t,x)$ the Green function of the dissipative transport operator given by (3.82). Then, we have the decomposition*

$$(3.98) \quad \Gamma(t,x) = K(t,x)\chi\{\underline{\lambda}t \leq x \leq \bar{\lambda}t, t \geq 1\} + \mathcal{K}(t,x) + R(t,x)\chi\{\underline{\lambda}t \leq x \leq \bar{\lambda}t\},$$

*where $R(t,x)$ can be written as*

$$(3.99) \quad \begin{aligned} R(t,x) = \sum_j \frac{e^{-(x-\lambda_j^1 t)^2/ct}}{1+t}\Big(&R_0(O(1))L_0 + R_0(O(1)(1+t)^{-1/2})L_- \\ &+ R_-(O(1)(1+t)^{-1/2})L_0 + R_-(O(1)(1+t)^{-1})L_-\Big), \end{aligned}$$

*for some constant $c$. Here $O(1)$ denotes a generic bounded matrix, $\lambda_j^1$ are the eigenvalues of the symmetric block $A_{11}$ of $A$, used in expansion (3.15), and the projectors are given by (3.19) and (3.20).*

## 4. The multi dimensional Green function

In this section we prove an analogous theorem for multi dimensional systems. Since, in general, the form of the Green function is not explicit, we have to relay directly on the Fourier coordinates. Thus the separation of the Green kernel into various part is done at the level of solution operator $\Gamma(t)$ acting on $L^2(R^m, R^n)$, or $L^1 \cap L^2(R^m, R^n)$. In the following we will consider the last space, even if one can study the equation for initial data only in $L^2$. Our aim is in fact to obtain decay estimates.

4.1. **General setting and first estimates.** We consider the Cauchy problem for the linear relaxation system in the Conservative-Dissipative form

$$w_t + \sum_{\alpha=1}^{m} A_\alpha w_{x_\alpha} = Bw, \qquad w \in \mathbb{R}^{n_1+n_2}, \tag{4.1}$$

$$w(0, \cdot) = w^0. \tag{4.2}$$

We assume that $A_\alpha$, $\alpha = 1, \ldots, m$, are symmetric matrices and that we have, as in (3.2),

$$B = \begin{bmatrix} 0 & 0 \\ 0 & D \end{bmatrix},$$

where $D$ is a negative definite matrix $\in \mathbb{R}^{n_2 \times n_2}$. So we have (H1).

Set, for $\xi \in \mathbb{R}^m$,

$$A(\xi) := \sum_{\alpha=1}^{m} \xi_\alpha A_\alpha, \quad E(i\xi) = B - iA(\xi). \tag{4.3}$$

We assume also that we have (H2), and we recall that, from Lemma 2.4, there exists a $c > 0$ such that if $\lambda(i\xi)$ is an eigenvalue of $E(i\xi)$, with $\xi \in \mathbb{R}^m$, then

$$\Re(\lambda(i\xi)) \leq -c \frac{|\xi|^2}{1 + |\xi|^2}. \tag{4.4}$$

Let us introduce the polar coordinates in $\mathbb{R}^m$

$$\xi = \rho\zeta, \ \rho = |\xi|, \ \zeta \in S^{m-1}$$

and set

$$E(i\rho, \zeta) = E(i\rho\zeta).$$

More generally, in $\mathbb{C} \otimes S^{m-1}$,

$$E(z, \zeta) = E(z\zeta) = B - zA(\zeta). \tag{4.5}$$

Moreover, since $S^{m-1}$ is compact, then when $E(z, \zeta)$ is considered in $\mathbb{C} \otimes S^{m-1}$, the points $z = 0$, $z = \infty$ are uniformly isolated exceptional point for all $\zeta$, while in general there are a finite number of exceptional curves for $0 < |z| < \infty$. Thus we can expand $E(z, \zeta)$ near $z = 0$ and $z = \infty$ as in the one dimensional case.

As before, we want to study the Green kernel $\Gamma(t, x)$ of (4.1). We recall that the support of $\Gamma$ is contained in the wave cone of (4.1), so that, for $t \geq 0$, $\Gamma(t, \cdot)$ has compact support. The solution of the Cauchy problem (4.1)-(4.2) is given by

$$w(t, \cdot) = \Gamma(t, \cdot) * w^0 \tag{4.6}$$

and, using the Fourier transform, we have

$$\hat{w}(t, \xi) = \hat{\Gamma}(t, \xi)\hat{w}^0(\xi) = e^{E(i\xi)t}\hat{w}^0(\xi). \tag{4.7}$$

We now use (4.4) to obtain our first decay estimates. For $a > 0$, we have:

$$\left|\chi(|\xi| > a) e^{E(i\xi)t}\right| \leq Ce^{-c\frac{a^2}{1+a^2}t}, \tag{4.8}$$

$$\left|\chi(|\xi| \leq a) e^{E(i\xi)t}\right| \leq Ce^{-\frac{c}{1+a^2}|\xi|^2 t}. \tag{4.9}$$

We have the following natural decomposition

$$w(t, \cdot) = M_a(t)w^0 + \mathcal{M}_a(t)w^0,$$

with

$$\widehat{\mathcal{M}_a(t)w^0} = \chi(|\xi| > a) e^{E(i\xi)t}\hat{w}^0(\xi), \tag{4.10}$$

$$\widehat{M_a(t)w^0} = \chi(|\xi| \leq a) e^{E(i\xi)t}\hat{w}^0(\xi). \tag{4.11}$$

For the high frequencies we obtain

$$\|\mathcal{M}_a(t)w^0\|_{L^2} = C\|\chi(|\xi| > a)\,e^{E(i\xi)t}\hat{w}^0(\xi)\|_{L^2}$$

$$\leq Ce^{-c\frac{a^2}{1+a^2}t}\|w^0\|_{L^2}$$

and, more generally, for any derivative $D^\beta$ in the space variables:

(4.12) $$\|D^\beta \mathcal{M}_a(t)w^0\|_{L^2} \leq Ce^{-c\frac{a^2}{1+a^2}t}\|D^\beta w^0\|_{L^2}.$$

On the other hand, for the low frequencies, we have

$$\|D^\beta M_a(t)w^0\|_{L^\infty} \leq C \int_{S^{m-1}} \int_0^a e^{-\frac{c}{1+a^2}|\xi|^2 t} |\xi|^\beta |\hat{w}^0(\xi)| |\xi|^{m-1} d|\xi| d\zeta$$

$$\leq C(a,|\beta|)\min\left(1, t^{-\frac{m}{2}-\frac{|\beta|}{2}}\right)\|w^0\|_{L^1}$$

and

$$\|D^\beta M_a(t)w^0\|_{L^2} \leq C\left(\int_{S^{m-1}} \int_0^a e^{-\frac{c}{1+a^2}|\xi|^2 t} |\xi|^{2\beta} |\hat{w}^0(\xi)|^2 |\xi|^{m-1} d|\xi| d\zeta\right)^{\frac{1}{2}}$$

$$\leq C(a,|\beta|)\min\left(1, t^{-\frac{m}{4}-\frac{|\beta|}{2}}\right)\|w^0\|_{L^1}$$

More generally, for $\beta \in \mathbb{N}^m$ and $p \in [2,+\infty]$, we obtain the decay estimates:

(4.13) $$\|D^\beta M_a(t)w^0\|_{L^p} \leq C(a,|\beta|)\min\left(1, t^{-\frac{m}{2}(1-\frac{1}{p})-\frac{|\beta|}{2}}\right)\|w^0\|_{L^1}.$$

To obtain a more refined estimate, we have to use the Conservative-Dissipative form in (4.13), by expanding $E(i\xi)$ for the low frequencies.

4.2. **Low frequencies estimates.** We now study the expansion of $E(z,\zeta) = B - zA(\zeta)$ near $z = 0$. We can use the result of Section 3, noting that the matrix $A$ in (3.5) is simply replaced by $A(\zeta)$. We introduce the total projector $P(z,\zeta)$ corresponding to all the eigenvalues near 0, and $P_-(z,\zeta) = I - P(z,\zeta)$ is the projector corresponding to the whole family of the eigenvalues with strictly negative real part (see (3.17) and (3.18)). The principal part of $P(z,\zeta)$ is the projector $Q_0 = R_0 L_0$, the principal part of $P_-(z,\zeta)$ is the projector $Q_- = R_- L_-$, and the projectors $R_0, L_0, R_-, L_-$ are given by

(4.14) $$L_0 = R_0^T = \begin{bmatrix} I_{n_1} & 0 \end{bmatrix}, \quad L_- = R_-^T = \begin{bmatrix} 0 & I_{n_2} \end{bmatrix}.$$

As in Section 3.1, we can write the expansion of the eigenprojectors $L(z,\zeta), R(z,\zeta)$ corresponding to the vanishing eigenvalues. By formula (3.29) we obtain

(4.15)
$$\begin{aligned}
L(z,\zeta) &= L_0 + zL_0 A(\zeta) R_- D^{-1} L_- + O(z^2) \\
&= L_0 + z\begin{bmatrix} 0 & A_{12}(\zeta)D^{-1} \end{bmatrix} + O(z^2) \\
&= \begin{bmatrix} I_{n_1} & zA_{12}(\zeta)D^{-1} \end{bmatrix} + O(z^2).
\end{aligned}$$

(4.16)
$$\begin{aligned}
R(z,\zeta) &= R_0 + zR_- D^{-1} L_- A(\zeta) R_0 + O(z^2) \\
&= R_0 + z\begin{bmatrix} 0 \\ D^{-1} A_{21}(\zeta) \end{bmatrix} + O(z^2) \\
&= \begin{bmatrix} I_{n_1} \\ zD^{-1} A_{21}(\zeta) \end{bmatrix} + O(z^2).
\end{aligned}$$

Thus, as in (3.32), we obtain

$$F(z, \zeta) \doteq L(z, \zeta)E(z, \zeta)R(z, \zeta)$$
$$= -zA_{11}(\zeta) - z^2 A_{12}(\zeta)D^{-1}A_{21}(\zeta) + O(z^3)$$
(4.17)
$$= \overline{F}(z, \zeta) + O(z^3).$$

In the same way, using (3.39) and (3.40), we obtain

(4.18)
$$F_-(z) \doteq L_-(z, \zeta)E(z, \zeta)R_-(z, \zeta)$$
$$= D - zA_{22}(\zeta) + O(z^2).$$

Let us recall that near $z = 0$ we have

$$P(z, \zeta) = R(z, \zeta)L(z, \zeta), \ L(z, \zeta)R(z, \zeta) = I,$$
$$P_-(z, \zeta) = R_-(z, \zeta)L_-(z, \zeta), \ L_-(z, \zeta)R_-(z, \zeta) = I,$$

and

$$E(z, \zeta) = R(z, \zeta)F(z, \zeta)L(z, \zeta) + R_-(z, \zeta)F_-(z, \zeta)L_-(z, \zeta).$$

This yields

(4.19)
$$e^{E(i\xi)t} = R(i\xi)e^{F(i\xi)t}L(i\xi) + R_-(i\xi)e^{F_-(i\xi)t}L_-(i\xi).$$

Take now a constant $a$ small enough, such that we can use decomposition (4.19) in (4.11):

$$\widehat{M_a(t)w^0} = \chi(|\xi| \leq a)R(i\xi)e^{F(i\xi)t}L(i\xi) + \chi(|\xi| \leq a)R_-(i\xi)e^{F_-(i\xi)t}L_-(i\xi)$$
(4.20)
$$\doteq \widehat{M_{a,1}(t)w^0} + \widehat{M_{a,2}(t)w^0}.$$

**Lemma 4.1.** *Assume $a \ll 1$. There exist two constants $c, C > 0$, such that*

(4.21)
$$\left|\chi(|\xi| \leq a)e^{F_-(i\xi)t}\right| \leq Ce^{-ct},$$

(4.22)
$$\left|\chi(|\xi| \leq a)e^{F(i\xi)t}\right| \leq Ce^{-c|\xi|^2 t}.$$

*Proof.* The matrix $D$ is negative definite, thus (4.18) implies (4.21). To obtain (4.22), we use again the polar coordinates $\rho = |\xi|$, $\xi = \rho\zeta$, to write

$$F(i\xi) = -i\rho A_{11}(\zeta) + \rho^2 A_{12}(\zeta)D^{-1}A_{21}(\zeta) + O(\rho^3)$$
$$= \overline{F}(i\xi) + O(\rho^3).$$

The matrix $A_{11}(\zeta)$ is real symmetric, and, by Lemma 3.1, the matrix $A_{12}(\zeta)D^{-1}A_{21}(\zeta)$ is negative defined (uniformly in $S^{m-1}$). If $\mu$ is an eigenvalue of the matrix

$$-iA_{11}(\zeta) + \rho A_{12}(\zeta)D^{-1}A_{21}(\zeta),$$

there exists a constant $c > 0$ such that

$$\Re(\mu) \leq -c\rho, \qquad (0 \leq \rho \leq a).$$

Therefore, if $\mu(i\xi)$ is an eigenvalue of $F(i\xi)$, we have

(4.23)
$$\Re(\mu(i\xi)) \leq -c|\xi|^2, \qquad (0 \leq |\chi| \leq a),$$

which yields (4.22). Let us underline that these inequalities are a direct consequence of assumption (H2). □

Next, we fix $a > 0$, such that estimates (4.21) and (4.22) hold, and we make a decomposition of the Green operator:

(4.24)
$$\Gamma(t) = K(t) + \mathcal{K}(t),$$

where

$$K(t) \doteq M_{a,1}(t), \quad \mathcal{K}(t) \doteq M_{a,2}(t) + \mathcal{M}_a(t).$$

For every function $w^0 \in L^1 \cap L^2(\mathbb{R}^m, \mathbb{R}^n)$, the solution $w(t) = \Gamma(t)w^0$ of (4.1), (4.2) can be decomposed as

$$w(t) = \Gamma(t)w^0 = K(t)w^0 + \mathcal{K}(t)w^0. \tag{4.25}$$

Using (4.12) and (4.21), we can estimate the second term on the RHS: there exist two constants $c, C > 0$, such that, for all space derivative $D^\beta$, it holds

$$\|D^\beta \mathcal{K}(t)w^0\|_{L^2} \le C e^{-ct} \|D^\beta w^0\|_{L^2}. \tag{4.26}$$

We can now establish the decay properties of $K(t)$, using the Conservative-Dissipative form. Using (4.15), (4.16), (4.20), and (4.22), we have that there exist two constants $c, C > 0$, such that

$$\begin{aligned}
\left|L_0 \widehat{K(t)w^0}\right| &\le C e^{-c|\xi|^2 t}\left(|L_0 \hat{w}^0| + |\xi| |L_- \hat{w}^0|\right) \\
\left|L_- \widehat{K(t)w^0}\right| &\le C e^{-c|\xi|^2 t}\left(|\xi| |L_0 \hat{w}^0| + |\xi|^2 |L_- \hat{w}^0|\right).
\end{aligned} \tag{4.27}$$

Using (4.27) we obtain

$$\begin{aligned}
\|L_0 K(t)w^0\|_{L^2}^2 &\le C \int_{\mathbb{S}^{m-1}} \int_0^\infty e^{-2c|\xi|^2 t}\left(|L_0 \hat{w}^0(\xi)|^2 + |\xi|^2 |L_- \hat{w}^0(\xi)|^2\right)|\xi|^{m-1} d|\xi| d\zeta \\
&\le C \min\{1, t^{-m/2}\}\|L_0 \hat{w}^0\|_{L^\infty}^2 + C \min\{1, t^{-m/2-1}\}\|L_- \hat{w}^0\|_{L^\infty}^2 \\
&\le C \min\{1, t^{-m/2}\}\|L_0 w^0\|_{L^1}^2 + C \min\{1, t^{-m/2-1}\}\|L_- w^0\|_{L^1}^2,
\end{aligned}$$

$$\begin{aligned}
\|L_- K(t)w^0\|_{L^2}^2 &\le C \int_{\mathbb{S}^{m-1}} \int_0^\infty e^{-2c|\xi|^2 t}\left(|\xi|^2 |L_0 \hat{w}^0(\xi)|^2 + |\xi|^4 |L_- \hat{w}^0(\xi)|\right)|\xi|^{m-1} d|\xi| d\zeta \\
&\le C \min\{1, t^{-m/2-1}\}\|L_0 w^0\|_{L^1}^2 + C \min\{1, t^{-m/2-2}\}\|L_- w^0\|_{L^1}^2.
\end{aligned}$$

Similarly, it is easy to prove that for every multi index $\beta$ the coefficient $\xi^{2\beta}$ appears in the integrand, so that

$$\|L_0 D^\beta K(t)w^0\|_{L^2} \le C \min\{1, t^{-m/4-|\beta|/2}\}\|L_0 w^0\|_{L^1} + C \min\{1, t^{-m/4-1/2-|\beta|/2}\}\|L_- w^0\|_{L^1}.$$

$$\|L_- D^\beta K(t)w^0\|_{L^2} \le C \min\{1, t^{-m/4-1/2-|\beta|/2}\}\|L_0 w^0\|_{L^1} + C \min\{1, t^{-m/4-1-|\beta|/2}\}\|L_- w^0\|_{L^1}.$$

We can estimate also the decay in every $p \in [2, +\infty]$. We have that

$$\begin{aligned}
|D^\beta L_0 K(t)w^0| &\le C \int e^{-c|\xi|^2 t} |\xi|^\beta \left(\|L_0 \hat{w}^0\|_{L^\infty} + |\xi| \|L_- \hat{w}^0\|_{L^\infty}\right) d\xi \\
&\le C \min\{1, t^{-m/2-|\beta|/2}\}\|L_0 w^0\|_{L^1} + C \min\{1, t^{-m/2-1/2-|\beta|/2}\}\|L_- w^0\|_{L^1},
\end{aligned}$$

$$|D^\beta L_- K(t)w^0| \le C \min\{1, t^{-m/2-1/2-|\beta|/2}\}\|L_0 w^0\|_{L^1} + C \min\{1, t^{-m/2-1-|\beta|/2}\}\|L_- w^0\|_{L^1},$$

so that, if $\beta$ is a multi index, for $p \in [2, +\infty]$ we have also the "$K(t)$ estimates":

$$\begin{aligned}
\|L_0 D^\beta K(t)w^0\|_{L^p} &\le C(|\beta|) \min\left\{1, t^{-\frac{m}{2}(1-\frac{1}{p})-|\beta|/2}\right\}\|w_c^0\|_{L^1} \\
&\quad + C(|\beta|) \min\left\{1, t^{-\frac{m}{2}(1-\frac{1}{p})-1/2-|\beta|/2}\right\}\|w_d^0\|_{L^1},
\end{aligned} \tag{4.28}$$

$$\begin{aligned}
\|L_- D^\beta K(t)w^0\|_{L^p} &\le C(|\beta|) \min\left\{1, t^{-\frac{m}{2}(1-\frac{1}{p})-1/2-|\beta|/2}\right\}\|w_c^0\|_{L^1} \\
&\quad + C(|\beta|) \min\left\{1, t^{-\frac{m}{2}(1-\frac{1}{p})-1-|\beta|/2}\right\}\|w_d^0\|_{L^1}.
\end{aligned} \tag{4.29}$$

4.3. **Decay estimates.** We thus collect the results in the following theorem.

**Theorem 4.2.** *Consider the linear PDE in the Conservative-Dissipative form*

$$w_t + \sum_{\alpha=1}^{m} A_\alpha w_{x_\alpha} = Bw, \tag{4.30}$$

*where $A_\alpha$, $B$ satisfy the assumption* **(SK)** *of Definition 2.3, and let $Q_0 = R_0 L_0$, $Q_- = I - Q_0 = R_- L_-$ be the eigenprojectors on the null space and the negative definite part of $B$. Then, for any function $w^0 \in L^1 \cap L^2(\mathbb{R}^m, \mathbb{R}^n)$, the solution $w(t) = \Gamma(t)w^0$ of (4.1), (4.2) can be decomposed as*

$$w(t) = \Gamma(t)w^0 = K(t)w^0 + \mathcal{K}(t)w^0, \tag{4.31}$$

*where the following estimates hold: for any multi index $\beta$ and for every $p \in [2, +\infty]$,*

$K(t)$ **estimates:**

$$\|L_0 D^\beta K(t)w^0\|_{L^p} \leq C(|\beta|) \min\left\{1, t^{-\frac{m}{2}(1-\frac{1}{p})-|\beta|/2}\right\} \|L_0 w^0\|_{L^1}$$
$$+ C(|\beta|) \min\left\{1, t^{-\frac{m}{2}(1-\frac{1}{p})-1/2-|\beta|/2}\right\} \|L_- w^0\|_{L^1}, \tag{4.32}$$

$$\|L_- D^\beta K(t)w^0\|_{L^p} \leq C(|\beta|) \min\left\{1, t^{-\frac{m}{2}(1-\frac{1}{p})-1/2-|\beta|/2}\right\} \|L_0 w^0\|_{L^1}$$
$$+ C(|\beta|) \min\left\{1, t^{-\frac{m}{2}(1-\frac{1}{p})-1-|\beta|/2}\right\} \|L_- w^0\|_{L^1}; \tag{4.33}$$

$\mathcal{K}(t)$ **estimates:**

$$\|D^\beta \mathcal{K}(t)w^0\|_{L^2} \leq C e^{-ct} \|D^\beta w^0\|_{L^2}. \tag{4.34}$$

*Remark* 4.3. Let us notice the relation among the Green kernel for (4.1) and the parabolic $n_1 \times n_1$ system in $m$ dimensions

$$w_t + \sum_{\alpha=1}^{m} A_{\alpha,11} w_{x_\alpha} = -\sum_{\alpha,\beta=1}^{m} A_{\alpha,12} D^{-1} A_{\beta,21} w_{x_\alpha x_\beta}, \qquad w \in \mathbb{R}^{n_1}. \tag{4.35}$$

This relation will be exploited better in the Sections 5.4, 5.5. Here we want to prove that the above system satisfies the following assumptions:

(1) there exists a unitary matrices $C(\zeta)$, $\zeta \in \mathcal{S}_{n-1}$, such that

$$C^T(\zeta)\left(\sum_{\alpha,\beta=1}^{m} \zeta_\alpha \zeta_\beta A_{\alpha,12} D^{-1} A_{\beta,21}\right) C(\zeta) = \begin{bmatrix} 0 & 0 \\ 0 & \hat{D}(\zeta) \end{bmatrix}, \tag{4.36}$$

with $\hat{D}(\zeta)$ negative definite (also its dimension depends on $\zeta$, in general);
(2) any eigenvector of $\sum_\alpha \zeta_\alpha A_{\alpha,11}$ is not in the null eigenspace of

$$\sum_{\alpha,\beta=1}^{m} \zeta_\alpha \zeta_\beta A_{\alpha,12} D^{-1} A_{\beta,21}.$$

It is easy to verify that the above assumptions correspond to Shizuta-Kawashima condition along any direction $\zeta$ for the parabolic system (4.35). To prove (1), let $C(\zeta) \in \mathbb{R}^{n_1 \times n_1}$ be the change of coordinates so that

$$\left(\sum_{\alpha=1}^{m} \zeta_\alpha A_{\alpha,21}\right) C(\zeta) = \begin{bmatrix} 0 & K(\zeta) \end{bmatrix},$$

with $\operatorname{Ker}(K(\zeta)) = \{0\}$. Since $A_{\alpha,12} = A_{\alpha,12}^T$, we thus obtain

$$C^T(\zeta)\left(\sum_{\alpha=1}^{m} \zeta_\alpha A_{\alpha,21}\right)^T D^{-1}\left(\sum_{\alpha=1}^{m} \zeta_\alpha A_{\alpha,21}\right) C(\zeta) = \begin{bmatrix} 0 \\ K^T(\zeta) \end{bmatrix} D^{-1} \begin{bmatrix} 0 & K(\zeta) \end{bmatrix}$$
$$= \begin{bmatrix} 0 & 0 \\ 0 & K^T(\zeta) D^{-1} K(\zeta) \end{bmatrix}.$$

The matrix $K^T(\zeta) D^{-1} K(\zeta)$ is negative definite, since $D^{-1}$ is negative definite and $\operatorname{Ker}(K(\zeta)) = \{0\}$.

Assume now that $v(\zeta)$, eigenvector of $\sum_\alpha \zeta_\alpha A_{\alpha,11}$, is in the null space of the viscosity matrix of (4.35). The it follows that is in $\text{Ker}(\sum_\alpha \zeta_\alpha A_{\alpha,21})$, so that the vector $R_0 v$ is an eigenvector of $\sum_\alpha \zeta_\alpha A_\alpha$. But this contradicts our assumptions, because is in the null space of $B$. For a related discussion about this remark, in a slightly different framework, see [38].

*Remark* 4.4. We note here that we cannot expect any estimate of the form

$$\|w(t)\|_{L^1} \leq C\|w^0\|_{L^1},$$

because for large $t$ the function $L_0 w$ behaves like the solution to

$$w_{c,t} + \sum_{\alpha=1}^m L_0 A_\alpha R_0 w_{c,x_\alpha} = 0,$$

and it is knows from [4] that this estimate is not true in general. The $L^\infty$ estimate depends strongly on the presence of a uniform parabolic operator, so that it is lost in the hyperbolic limit.

*Remark* 4.5. We note here that by means of the explicit form of the kernel in the one dimensional case it follows that the decay estimates holds for $p \in [1, +\infty]$, with the same decay rate.

*Remark* 4.6. Since in general we are not able to give the explicit form of the kernel part $K(t)$, one may suspect that even if the kernel $K(t)R_- L_-$ has the same decay estimates of a derivative of the heat kernel, it is not a derivative of a heat like kernel. This is striking different from the one dimensional case.

However, a simple observation shows that the function $L_0 \Gamma(t, x) R^-$ is actually a derivative. Note that in one space dimension we only proved that its principal part $L_0 K(t, x) R_-$ is an $x$-derivative. Thus we are obtaining a new result also in one space dimension.

By replacing $w_d$ with $\tilde{w}_d + e^{Dt} w_d(t = 0)$, we obtain that the equations for $(w_c, \tilde{w}_d)$ are

$$(4.37) \quad \begin{cases} w_{c,t} + \sum_\alpha A_{\alpha,11} w_{c,x_\alpha} + \sum_\alpha A_{\alpha,12} \tilde{w}_{d,x_\alpha} = -\sum_\alpha A_{\alpha,12} e^{Dt} w_{d,x_\alpha}(t=0) \\ \tilde{w}_{d,t} + \sum_\alpha A_{\alpha,21} w_{c,x_\alpha} + A_{\alpha,22} \tilde{w}_{d,x_\alpha} = D\tilde{w}_d - \sum_\alpha A_{\alpha,22} e^{Dt} w_{d,x_\alpha}(t=0) \end{cases}$$

with initial data $(w_c(t=0), 0)$. The solution can be written by Duhamel formula as

$$\begin{pmatrix} w_c \\ \tilde{w}_d \end{pmatrix} = \Gamma(t) * \begin{pmatrix} w_c(t=0) \\ 0 \end{pmatrix} - \sum_\alpha \int_0^t \Gamma(t-s) * \begin{pmatrix} A_{\alpha,12} e^{Ds} w_{d,x_\alpha}(t=0) \\ A_{\alpha,22} e^{Ds} w_{d,x_\alpha}(t=0) \end{pmatrix} ds$$

$$(4.38) \quad = \Gamma(t) * \begin{pmatrix} w_c(t=0) \\ 0 \end{pmatrix} - \sum_\alpha \partial_{x_\alpha} \int_0^t \Gamma(t-s) * \begin{pmatrix} A_{\alpha,12} e^{Ds} w_d(t=0) \\ A_{\alpha,22} e^{Ds} w_d(t=0) \end{pmatrix} ds.$$

In particular, one sees that the

$$(4.39) \quad \Gamma(t,x) R_- = \sum_\alpha \partial_{x_\alpha} \left( \int_0^t \Gamma(t-s,x) A_\alpha R_- e^{Ds} ds \right) + \begin{pmatrix} 0 \\ e^{Dt} \delta(x) \end{pmatrix}.$$

This remark is useful to deal with the case $m = 2$ in Section 5.4.

*Example* 4.7. **Rotationally invariant systems.** If we assume that system (4.30) is invariant for rotations, it is possible to give a more precise expansion of the parabolic part $K(t)$ of the kernel $\Gamma(t)$. Consider for example the linearized isentropic Euler equations with damping, which can be written as

$$(4.40) \quad \begin{cases} \rho_t + \text{div} v = 0 \\ v_t + \nabla \rho = -v \end{cases}$$

To fix the ideas, take $m = 3$, $n = 4 = n_1 + n_2 = 1 + 3$. Clearly the system is already in the Conservative-Dissipative form and the condition **(SK)** is satisfied. In this case one can decompose $K(t, x)$ as

$$(4.41) \quad K(t,x) = \begin{bmatrix} G(t,x) & (\nabla G(t,x))^T \\ \nabla G(t,x) & \nabla^2 G(t,x) \end{bmatrix} + R_1(t,x),$$

where $G(t, x)$ is the heat kernel for $u_t = \Delta u$, and the rest term $R_1(t, x)$ satisfies the bound

$$(4.42) \quad R_1(t,x) = \frac{e^{-c|x|^2/t}}{(1+t)^2} \begin{bmatrix} O(1) & O(1)(1+t)^{-1/2} \\ O(1)(1+t)^{-1/2} & O(1)(1+t)^{-1} \end{bmatrix}.$$

In particular the principal part of $\Gamma(t)$ is given by the heat kernel $G(t,x)$.

A more interesting example is the system

$$(4.43) \quad \begin{cases} \rho_t + \text{div}\, v = 0, \\ v_t + \nabla \rho + \text{div}\, R = 0, \\ R_t + \nabla v = -R, \end{cases}$$

where $\rho \in \mathbb{R}$, $v \in \mathbb{R}^3$, and $R \in \mathbb{R}^9$. In this case, thanks to the invariance for rotations of the Green kernel, we can use the one dimensional decomposition (3.57) to find that the main smooth part $K_{00}(t)$ of the Green kernel $\Gamma(t)$ is given by

$$(4.44) \quad K_{00}(t) = \begin{bmatrix} 0 & 0 \\ 0 & G(t,x)\mathcal{P} \end{bmatrix} + \begin{bmatrix} (W_{00} * G)(t,x) & (W_{01} * G)(t,x) \\ (W_{10} * G)(t,x) & (W_{11} * G)(t,x) \end{bmatrix} + R_1(t,x).$$

Here $\mathcal{P} : (L^2(\mathbb{R}^3))^3 \mapsto (L^2(\mathbb{R}^3))^3$ is the orthogonal projection of $L^2$ vector fields on the subspace of divergence free vector fields. $\mathcal{P}v$ is characterized by

$$\mathcal{P}v \in (L^2(\mathbb{R}^3))^3, \quad \text{div}\,\mathcal{P}v = 0, \quad \text{curl}(v - \mathcal{P}v) = 0,$$

and so we have that

$$v - \mathcal{P}v = \nabla \psi, \quad \text{with } \Delta \psi = \text{div}\, v.$$

This yields

$$(4.45) \quad \mathcal{P}v = v - \nabla(\Delta^{-1} \text{div}\, v).$$

In fact, in Fourier coordinates, we have

$$(4.46) \quad \widehat{\mathcal{P}v}(\xi) = \hat{v}(\xi) - |\xi|^{-2}(\xi \cdot \hat{v}(\xi))\xi = \hat{v}(\xi) - |\xi|^{-2}\xi\xi^T \cdot \hat{v}(\xi).$$

The matrix valued function

$$(4.47) \quad W(t,x) = W_1(t,x) + W_2(t,x) = \begin{bmatrix} W_{00}(t,x) & W_{01}(t,x) \\ W_{10}(t,x) & W_{11}(t,x) \end{bmatrix} + \begin{bmatrix} 0 & 0 \\ 0 & \delta(x)\mathcal{P} \end{bmatrix}$$

is the matrix valued Green function of the system

$$\begin{cases} \rho_t + \text{div}\, v = 0 \\ v_t + \nabla \rho = 0 \end{cases}$$

and it can be written by means of the fundamental solution to the wave equation $u_{tt} = \Delta u$. In fact, $W_{00}$ is the solution of $u_{tt} = \Delta u$ with initial data $u = \delta(x)$, $u_t = 0$, and

$$(4.48) \quad W_1 = \begin{bmatrix} W_{00} & \nabla^T \partial_t (-\Delta)^{-1} W_{00} \\ \nabla \partial_t (-\Delta)^{-1} W_{00} & -\nabla^2 (-\Delta)^{-1} W_{00} \end{bmatrix}.$$

In particular one can check that $W_2$ corresponds to incompressible vector fields, while $W_1$ corresponds to curl free vector fields. Finally the rest $R_1(t,x)$ satisfies

$$|R_1(t,x)| \leq (1+t)^{-1/2} |G(t,x)| + (1+t)^{-1/2} |(W * G)(t,x)|.$$

From (4.44), one sees that the asymptotic behaviour of $\Gamma(t)$ is a function $(0, v_0)$, with $v_0$ divergence free vector field, which remains close to the origin, and a function $(\rho, v_1)$, with $v_1$ curl free, which diffuses around the sound cone $\{|x| = t\}$. Due to the finite speed of propagation of (4.43), we can restrict $K_{00}(t)$ to the light cone $\{|x| \leq \sqrt{2}t\}$. Finally, let us notice that the main part of the kernel $K_{00}$ is the Green function of the fully parabolic system

$$(4.49) \quad \begin{cases} \rho_t + \text{div}\, v = \Delta \rho, \\ v_t + \nabla \rho = \Delta v. \end{cases}$$

## 5. Decay estimates in the nonlinear case and more accurate asymptotic behavior

In this section we study the time decay properties of the global smooth solutions to a nonlinear entropy strictly dissipative relaxation system in conservative-dissipative form. We shall prove that the conservative variables $u_c = L_0 u$ decays as the heat kernel and derivatives, while the dissipative variable $u_d = L_- u$ decays faster. Following (3.51), we set

$$K(t) = \begin{bmatrix} L_0 K(t) R_0 & L_0 K(t) R_- \\ L_- K(t) R_0 & L_- K(t) R_- \end{bmatrix} = \begin{bmatrix} K_{00}(t) & K_{0-}(t) \\ K_{-0}(t) & K_{--}(t) \end{bmatrix}.$$

Moreover we shall prove that $u_c(t)$ approaches the conservative part $K_{00}(t) L_0 u(0)$ of the linear solution $\Gamma(t) u(0)$ faster that the decay of the heat kernel for $m \geq 2$. In one dimension we shall show that $u_c(t)$ converges to the solution of a parabolic equation with quadratic nonlinearity, in the spirit of Chapman-Enskog expansion.

### 5.1. Decay estimates in $L^p$.

We now prove the decay estimates in $L^p(\mathbb{R}^m; \mathbb{R}^n)$, $p \in [2, +\infty]$, for the solution $u$ with initial data in $L^1 \cap H^s$, with $s$ sufficiently large, for the non linear equation

$$(5.1) \qquad u_t + \sum_{\alpha=1}^m (f_\alpha(u))_{x_\alpha} = g(u),$$

with $f_\alpha(0) = g(0) = 0$ and initial condition

$$(5.2) \qquad u(x, 0) = u_0(x).$$

We shall assume that the system (5.1) is strictly entropy dissipative and condition **(SK)** is satisfied. Under the assumptions of Theorem 2.5, we consider the global solution $u$ of (5.1)-(5.2), with

$$(5.3) \qquad u \in C^0\left([0, +\infty); H^s(\mathbb{R}^m)\right) \cap C^1\left([0, +\infty); H^{s-1}(\mathbb{R}^m)\right),$$

and we can assume that there exists $\delta_0 > 0$ such that, for $(t, x) \in [0, +\infty) \times \mathbb{R}^m$,

$$(5.4) \qquad |u(t, x)| \leq \delta_0.$$

Consider now the associated linearized problem

$$(5.5) \qquad u_t + \sum_{\alpha=1}^m Df_\alpha(0) u_{x_\alpha} = Dg(0) u.$$

Thanks to Proposition 2.7 and Remark 2.8, we can assume, without loss of generality, that this system is in the Conservative-Dissipative form (C-D form) of Definition 2.6. Therefore, thanks to Theorem 4.2, the associated Green function can be decomposed as

$$(5.6) \qquad \Gamma(t) = K(t) + \mathcal{K}(t),$$

and the estimates (4.32), (4.33), and (4.34) hold true.

We can write the solution $u$ to (5.1) by Duhamel's formula as

$$u(t) = \Gamma(t) u(0) + \sum_{\alpha=1}^m \int_0^t D_{x_\alpha} \Gamma(t-s) \bigl(Df_\alpha(0) u(s) - f_\alpha(u(s))\bigr) ds$$

$$(5.7) \qquad\qquad + \int_0^t \Gamma(t-s)(g(u(s)) - Dg(0) u(s)) ds.$$

*Remark* 5.1. We now observe that the only terms acting on $g(u) - Dg(0)u$ are the exponential decaying terms $\mathcal{K}(t)$, and the terms $K_{0-}(t)$ and $K_{--}(t)$. Thus when projecting on the conservative variables, the term with the slowest decay is $K_{0-}(t)$, while when projected on the dissipative variables the leading term is $K_{--}(t)$. A similar observation can be made for the product $\Gamma_{x_\alpha}(t-s)(Df_\alpha(0) u(s) - f_\alpha(u(s)))$, where the principal part for $u_c$ is $K_{00}(t)$ while for $u_d(t)$ is $K_{-0}(t)$.

We first prove that the solution decays as the heat kernels and derivatives, with no distinction among the conservative part $u_c = L_0 u$ and the dissipative one $u_d = L_- u$. We next prove that the dissipative part decays faster, as a derivative of the conservative one. Finally we shall study the decay of the time derivative.

For the $\beta$ derivative, we shall use the formula

$$D^\beta u(t) = D^\beta K(t)u(0) + \mathcal{K}(t)D^\beta u(0)$$
$$+ \sum_{\alpha=1}^{m} \int_0^{t/2} D^\beta D_{x_\alpha} K(t-s)\Big(Df_\alpha(0)u(s) - f_\alpha(u(s))\Big)ds$$
$$+ \sum_{\alpha=1}^{m} \int_{t/2}^{t} D_{x_\alpha} K(t-s)D^\beta\Big(Df_\alpha(0)u(s) - f_\alpha(u(s))\Big)ds$$
$$+ \int_0^{t/2} D^\beta K(t-s)R_- L_-\Big(g(u(s)) - Dg(0)u(s)\Big)ds$$
$$+ \int_{t/2}^{t} K(t-s)R_- D^\beta L_-\Big(g(u(s)) - Dg(0)u(s)\Big)ds$$
(5.8)
$$+ \sum_{\alpha=1}^{m} \int_0^{t} \mathcal{K}(t-s)D^\beta\Big(D_{x_\alpha}\big(Df_\alpha(0)u(s) - f_\alpha(u(s))\big) + \big(g(u(s)) - Dg(0)u(s)\big)\Big)ds,$$

observing that for $\beta = 0$ we do not need to split the integral in time.

Now we recall some well-known inequalities in the Sobolev spaces $H^s(\mathbb{R}^m)$. The basic remark is the fact that for $s > \frac{m}{2}$, $H^s(\mathbb{R}^m)$ is a Banach algebra, i.e., for any $u, v$ in $H^s(\mathbb{R}^m)$, we have:

(5.9) $$\|uv\|_{H^s} \leq C\|u\|_{H^s}\|v\|_{H^s}.$$

More generally, we need some Moser-type calculus inequalities, see for instance [23, 35].

For every $u, v$ in $H^s(\mathbb{R}^m) \cap L^\infty(\mathbb{R}^m)$ ($s \geq 0$), and $|\beta| \leq s$

(5.10) $$\|D^\beta(uv)\|_{L^2} \leq C\Big(\|u\|_{L^\infty}\|D^\beta v\|_{L^2} + \|v\|_{L^\infty}\|D^\beta u\|_{L^2}\Big),$$

so that, if $u, v \in H^{s+|\beta|}(\mathbb{R}^m)$,

(5.11) $$\|D^\beta(uv)\|_{H^s} \leq C\Big(\|u\|_{L^\infty}\|D^\beta v\|_{H^s} + \|v\|_{L^\infty}\|D^\beta u\|_{H^s}\Big).$$

For every smooth function $h : \mathbb{R} \mapsto \mathbb{R}$, every $u \in H^s(\mathbb{R}^m) \cap L^\infty(\mathbb{R}^m)$ ($s \geq 1$), which verifies inequality (5.4), and $\beta \neq 0$, $|\beta| \leq s$, we have

(5.12) $$\|D^\beta h(u)\|_{L^2} \leq C_\beta \|h'\|_{C^{|\beta|-1}(|u|\leq \delta_0)} \|u\|_{L^\infty}^{|\beta|-1} \|D^\beta u\|_{L^2}.$$

Moreover, if $h(0) = 0$, we have

(5.13) $$\|h(u)\|_{H^s} \leq C\Big(\delta_0, \|h\|_{C^s(|u|\leq \delta_0)}\Big)(1 + \|u\|_{H^s}).$$

As in [6], we use the following crude, but useful Lemma:

**Lemma 5.2.** *For any $\gamma, \delta \geq 0$, $t \geq 2$*

(5.14) $$\varphi \doteq \min\{\gamma, \delta, \gamma + \delta - 1\},$$

*then it holds*

(5.15) $$\int_0^t \min\{1, (t-s)^{-\gamma}\} \min\{1, s^{-\delta}\} ds \leq C \cdot \begin{cases} t^{-\varphi} & \gamma, \delta \neq 1 \\ t^{-\varphi}(1 + \ln t) & \gamma \leq 1, \delta = 1 \text{ or } \gamma = 1, \delta \leq 1 \\ t^{-1} & \gamma > 1, \delta = 1 \text{ or } \gamma = 1, \delta > 1 \end{cases}$$

(5.16) $$\int_0^t \min\{1, s^{-\delta}\} = C \cdot \begin{cases} 1 & \delta > 1 \\ \ln t & \delta = 1 \\ t^{1-\delta} & 0 \leq \delta < 1 \end{cases}$$

(5.17) $$\int_0^t e^{-c(t-s)} \min\{1, s^{-\gamma}\} \leq C \min\{1, s^{-\gamma}\}, \quad \gamma \geq 0,$$

5.1.1. *$H^s$ estimates.* Let us set now

$$E_s = \max\{\|u(0)\|_{L^1}, \|u(0)\|_{H^s}\}, \tag{5.18}$$

and

$$M_0(t) = \sup_{0 \leq \tau \leq t} \left\{\max\{1, \tau^{m/4}\} \|u(\tau)\|_{H^s}\right\}, \tag{5.19}$$

We are going to estimate the $H^s$ norm of $u(t)$ in (5.7).

For shortness we denote (the products below should be intended as tensor products)

$$f_\alpha(u) - Df_\alpha(0)u = u^2 h_\alpha(u), \qquad g(u) - Dg(0)u = u^2 h(u). \tag{5.20}$$

Using Theorem 4.2 and recalling Remark 5.1, we have the estimate

$$\|u(t)\|_{H^s} \leq \quad C \min\{1, t^{-m/4}\} \|u(0)\|_{L^1} + Ce^{-ct} \|u(0)\|_{H^s}$$

$$+ C \int_0^t \min\{1, (t-s)^{-m/4-1/2}\} \left(\sum_{\alpha=1}^m \|u^2 h_\alpha(u(s))\|_{L^1} + \|u^2 h(u(s))\|_{L^1}\right) ds \tag{5.21}$$

$$+ C \int_0^t e^{-c(t-s)} \left(\sum_{\alpha=1}^m \|D_{x_\alpha}(u^2 h_\alpha(u(s)))\|_{H^s} + \|u^2 h(u(s))\|_{H^s}\right) ds.$$

It is obvious that

$$C \min\{1, t^{-m/4}\} \|u(0)\|_{L^1} + Ce^{-ct} \|u(0)\|_{H^s} \leq C \min\{1, t^{-m/4}\} E_s.$$

For the third term, using (5.4), we have:

$$\sum_{\alpha=1}^m \|u^2 h_\alpha(u(s))\|_{L^1} + \|u^2 h(u(s))\|_{L^1}$$

$$\leq C \left(\sum_{\alpha=1}^m \|h_\alpha\|_{L^\infty(|u| \leq \delta_0)} + \|h\|_{L^\infty(|u| \leq \delta_0)}\right) \|u(s)\|_{L^2}^2 \tag{5.22}$$

$$\leq C \min\{1, s^{-m/2}\} (M_0(s))^2$$

Let us consider now the fourth term in (5.21). We have to estimate $\|D^\beta(u^2 h(u))\|_{H^s}$ for $\beta = 0$ and $|\beta| = 1$. More generally we have the following result.

**Lemma 5.3.** *Fix $s > \frac{m}{2}$ and $\beta \in \mathbb{N}^m$, and let $u \in H^r$ which verifies inequality (5.4), with $r \geq s + |\beta|$. Then we have*

$$\|D^\beta(u^2 h(u))\|_{H^s} \leq C\left(\delta_0, \|u\|_{H^s}, \|h\|_{C^{s+|\beta|}(|u| \leq \delta_0)}\right) \|u\|_{L^\infty} \|D^\beta u\|_{H^s}. \tag{5.23}$$

*Proof.* First we consider the case $\beta = 0$. We have

$$\|u^2 h(u)\|_{H^s} \leq \|u^2 h(0)\|_{H^s} + \|u^2(h(u) - h(0))\|_{H^s}.$$

Using (5.9) and (5.11) yields

$$\|u^2 h(u)\|_{H^s} \leq C(\|u\|_{L^\infty} \|u\|_{H^s} + \|u\|_{L^\infty} \|h(u) - h(0)\|_{L^\infty} \|u\|_{H^s}$$

$$+ \|u\|_{L^\infty} \|u\|_{H^s} \|h(u) - h(0)\|_{H^s}).$$

So, by (5.13), we obtain (5.23) for $\beta = 0$. For $|\beta| \geq 1$, using twice (5.11) yields

$$\|D^\beta(u^2 h(u))\|_{H^s} \leq C\|u\|_{L^\infty}\left(\|u\|_{L^\infty} \|D^\beta h(u)\|_{H^s} + \|h(u)\|_{L^\infty} \|D^\beta u\|_{H^s}\right).$$

Since $|\beta| \geq 1$, we can use (5.12) to obtain

$$\|D^\beta h(u)\|_{H^s} \leq C(\delta_0) \|h'\|_{C^{|\beta|+s-1}(|u| \leq \delta_0)} \|D^\beta u\|_{H^s},$$

which implies (5.23) with $C = C\left(\delta_0, \|h\|_{C^{s+|\beta|}(|u| \leq \delta_0)}\right)$. □

Now, inequality (5.23) gives

$$\sum_{\alpha=1}^{m} \|D_{x_\alpha}(u^2 h_\alpha(u))\|_{H^s} + \|u^2 h(u)\|_{H^s} \leq C\|u\|_{L^\infty}\|u\|_{H^{s+1}}. \tag{5.24}$$

Collecting inequalities (5.22) and (5.24) in (5.21), we obtain

$$\|u(t)\|_{H^s} \leq C \min\{1, t^{-m/4}\} E_s + C(M_0(t))^2 \int_0^t \min\{1, (t-s)^{-m/4-1/2}\} \min\{1, s^{-m/2}\} ds$$
$$+ CM_0(t) E_{s+1} \int_0^t e^{-c(t-s)} \min\{1, s^{-m/4}\} ds. \tag{5.25}$$

By means of Lemma 5.2, for $m = 1$ the first integral decays as $t^{-1/4}$, while for $m = 2$ decays as $t^{-1} \ln t < t^{-1/2}$, and for $m \geq 3$ as $t^{-m/4-1/2} < t^{-m/4}$. So, we obtain in (5.25)

$$M_0(t) \leq C\left(E_s + (M_0(t))^2 + E_{s+1} M_0(t)\right). \tag{5.26}$$

Then, if $E_{s+1}$ is small enough, we have the bound

$$M_0(t) \leq CE_s, \tag{5.27}$$

which implies

$$\|u(t)\|_{H^s} \leq C \min\{1, t^{-m/4}\} E_s. \tag{5.28}$$

5.1.2. $L^\infty$ *estimates.* We now estimate the $L^\infty$ norm of the solution. Set as before

$$N_0(t) = \sup_{0 \leq \tau \leq t} \left\{ \max\{1, \tau^{m/2}\} \|u(\tau)\|_{L^\infty} \right\}.$$

From (5.8) and Theorem 4.2, we have, using (5.22) and (5.27),

$$\|u(t)\|_{L^\infty} \leq C \min\{1, t^{-m/2}\} \|u(0)\|_{L^1} + Ce^{-ct} \|u(0)\|_{H^{[m/2]+1}}$$
$$+ CE_{[m/2]+1}^2 \int_0^t \min\{1, (t-s)^{-m/2-1/2}\} \min\{1, s^{-m/2}\} ds \tag{5.29}$$
$$+ CN_0(t) E_{[m/2]+1} \int_0^t e^{-c(t-s)} \min\{1, s^{-m/2}\} ds.$$

Using Lemma 5.2, for $m \geq 2$, we obtain

$$N_0(t) \leq C(E_{[m/2]+1} + E_{[m/2]+1} N_0(t) + E_{[m/2]+1}^2), \tag{5.30}$$

which implies

$$\|u(t)\|_{L^\infty} \leq C \min\{1, t^{-m/2}\} E_{[m/2]+1}, \tag{5.31}$$

if $E_{[m/2]+1}$ is small enough.

For $m = 1$, using the decomposition of Theorem 3.7 and the estimate (5.27), we can estimate the first integral as

$$\int_0^t \left(\|K_x(t)\|_{L^2} + \|R_x(t)\|_{L^2}\right) \|u^2\|_{L^2} \leq CN_0(t) E_1 \int_0^t \min\{1, (t-s)^{-3/4}\} \min\{1, s^{-3/4}\} ds$$
$$\leq C \min\{1, t^{-1/2}\} N_0(t) E_1,$$

and for $E_1$ is uniformly small, the inequality (5.31) holds also for $m = 1$.

5.1.3. *$D^\beta u$ estimates.* We set now for $b \in \mathbb{N}$

$$
(5.32) \qquad M_b(t) = \sum_{0 \leq |\beta| \leq b} \sup_{0 \leq \tau \leq t} \left\{ \max\left\{1, \tau^{m/4+|\beta|/2}\right\} \|D^\beta u(\tau)\|_{H^s} \right\},
$$

and we assume that for $0 \leq b < \bar{b}$ the following estimates hold:

$$
(5.33) \qquad M_b(t) \leq CE_{b+s} \ll 1.
$$

Now we estimate the $|\beta| = \bar{b}$ derivative using Theorem 4.2 with the decomposition (5.8) and Remark 5.1. We have

$$
\begin{aligned}
\|D^\beta u(t)\|_{H^s} \leq{}& C \min\{1, t^{-m/4-|\beta|/2}\} \|u(0)\|_{L^1} + Ce^{-ct}\|D^\beta u(0)\|_{H^s} \\
& + C \int_0^{t/2} \min\{1, (t-s)^{-m/4-1/2-|\beta|/2}\} \left( \sum_{\alpha=1}^m \|u^2 h_\alpha(u(s))\|_{L^1} + \|u^2 h(u(s))\|_{L^1} \right) ds \\
& + C \int_{t/2}^t \min\{1, (t-s)^{-m/4-1/2}\} \left( \sum_{\alpha=1}^m \|D^\beta D_{x_\alpha} u^2 h_\alpha(u(s))\|_{L^1} + \|D^\beta u^2 h(u(s))\|_{L^1} \right) ds \\
& + C \int_0^t e^{-c(t-s)} \left( \sum_{\alpha=1}^m \|D^\beta D_{x_\alpha} u^2 h_\alpha(u(s))\|_{H^s} + \|D^\beta u^2 h(u(s))\|_{H^s} \right) ds.
\end{aligned}
\qquad (5.34)
$$

For the first integral we obtain, using (5.22), (5.27), and Lemma 5.2:

$$
\begin{aligned}
\int_0^{t/2} & \min\{1, (t-s)^{-m/4-1/2-|\beta|/2}\} \left( \sum_{\alpha=1}^m \|u^2 h_\alpha(u(s))\|_{L^1} + \|u^2 h(u(s))\|_{L^1} \right) ds \\
& \leq CE_s^2 \int_0^{t/2} \min\{1, (t-s)^{-m/4-1/2-|\beta|/2}\} \min\{1, s^{-m/2}\} ds \\
& \leq C \min\{1, t^{-m/4-|\beta|/2}\} E_s^2.
\end{aligned}
\qquad (5.35)
$$

For the second integral we need to prove the following estimate

$$
(5.36) \qquad \left( \sum_{\alpha=1}^m \|D^\beta D_{x_\alpha} u^2 h_\alpha(u)\|_{L^1} + \|D^\beta u^2 h(u)\|_{L^1} \right) \leq C \min\{1, t^{-m/2-|\beta|/2}\} \left( E_{\bar{b}+s-1}^2 + E_s M_{\bar{b}} \right).
$$

We just consider the first term on the left. Using the chain rule, we have

$$
\begin{aligned}
D^\beta D_{x_\alpha} u^2 h_\alpha(u) &= D^\beta \left( (D_{x_\alpha} u) u \tilde{h}(u) \right) \\
&= \sum_{\alpha \leq \beta} \sum_{\gamma \leq \beta - \alpha} \sum_{\delta_1 + \cdots + \delta_k = \alpha} C(\alpha, \beta, \gamma, \delta_i)(D^{\beta-\alpha-\gamma} D_{x_\alpha} u)(D^\gamma u)(D^{\delta_1} u) \cdots (D^{\delta_k} u) \tilde{h}^{(k)}(u).
\end{aligned}
$$

Following the proof of Lemma 3.10 in [35], we obtain for the generic term

$$
\begin{aligned}
\|(D^{\beta-\alpha-\gamma} D_{x_\alpha} u)(D^\gamma u)(D^{\delta_1} u) \cdots (D^{\delta_k} u)\|_{L^1} &\leq \|(D^{\beta-\alpha-\gamma} D_{x_\alpha} u)\|_{L^2} \cdot \|(D^\gamma u)(D^{\delta_1} u) \cdots (D^{\delta_k} u)\|_{L^2} \\
&\leq C \|D^{\beta-\alpha-\gamma} u\|_{H^1} \|D^{|\gamma+\alpha|} u\|_{L^2}.
\end{aligned}
$$

Then, to obtain (5.36), we have to use (5.33) for the cases $\gamma + \alpha \neq 0$ and $\gamma + \alpha \neq \beta$.

For the third integral, the inequality (5.23) yields

$$
\begin{aligned}
\left( \sum_{\alpha=1}^m \|D^\beta D_{x_\alpha} u^2 h_\alpha(u(s))\|_{H^s} + \|D^\beta u^2 h(u(s))\|_{H^s} \right) &\leq C\|u\|_{L^\infty} \|D^\beta u\|_{H^{s+1}} \\
&\leq C\|u\|_{L^\infty} \|D^{\beta'} u\|_{H^{s+2}},
\end{aligned}
$$

where $|\beta'| = \bar{b} - 1$. We use the $L^\infty$-estimate (5.31) and the induction hypothesis (5.33) (replacing $s$ by $s+2$), to obtain

$$\text{(5.37)} \quad \left(\sum_{\alpha=1}^{m} \|D^\beta D_{x_\alpha} u^2 h_\alpha(u(s))\|_{H^s} + \|D^\beta u^2 h(u(s))\|_{H^s}\right) \leq C \min\{1, t^{-m/4-\frac{\bar{b}}{2}}\} E_{[m/2]+1} E_{\bar{b}+s+1}.$$

Substituting the above inequalities into (5.34) yields

$$
\begin{aligned}
\|D^\beta u(t)\|_{H^s} &\leq C \min\{1, t^{-m/4-|\beta|/2}\} \left(E_{\bar{b}+s} + E_s^2\right) \\
\text{(5.38)} \quad &+ C\left(E_{\bar{b}+s-1}^2 + E_s M_{\bar{b}}(t)\right) \int_{t/2}^{t} \min\{1, (t-s)^{-m/4-1/2}\} \min\{1, s^{-m/2-|\beta|/2}\} ds \\
&+ C E_{[m/2]+1} E_{\bar{b}+s+1} \int_0^t e^{-c(t-s)} \min\{1, s^{-m/4-\frac{\bar{b}}{2}}\} ds.
\end{aligned}
$$

Using again Lemma 5.2 we finally obtain

$$\text{(5.39)} \quad \|D^\beta u(t)\|_{H^s} \leq C \min\{1, t^{-m/4-|\beta|/2}\} \left(E_{\bar{b}+s} + E_s^2 + E_{\bar{b}+s-1}^2 + E_{[m/2]+1} E_{\bar{b}+s+1} + E_s M_{\bar{b}}(t)\right).$$

If $E_s$ is small enough, we have

$$M_{\bar{b}}(t) \leq C E_{\bar{b}+s},$$

with $C = C(E_{\bar{b}+s+1})$. So, we can conclude

$$\text{(5.40)} \quad \|D^\beta u(t)\|_{H^s} \leq C \min\{1, t^{-m/4-|\beta|/2}\} E_{|\beta|+s}.$$

To obtain the $L^\infty$ estimates we use the following inequalities from [35], Proposition 3.8:

$$\|D^\beta u(t)\|_{L^\infty} \leq C \|D^{[m/2]+1} D^\beta u(t)\|_{L^2}^{1/2} \|D^k D^\beta u(t)\|_{L^2}^{1/2},$$

with $k = [m/2] - 1$ if $m$ is even and $k = [m/2]$ if $m$ is odd. So, we just use (5.40) with $s = [m/2]+1$ to obtain

$$\text{(5.41)} \quad \|D^\beta u(t)\|_{L^\infty} \leq C \min\{1, t^{-m/2-|\beta|/2}\} E_{|\beta|+m+2}.$$

Actually, it is possible to show, by a direct calculation, that the following estimate holds:

$$\text{(5.42)} \quad \|D^\beta u\|_{L^\infty} \leq C \min\{1, t^{-m/2-|\beta|/2}\} E_{|\beta|+[m/2]+1}.$$

Therefore, we easily obtain the decay estimate in the $L^p$-spaces

$$\text{(5.43)} \quad \|D^\beta u(t)\|_{L^p} \leq C \min\{1, t^{-\frac{m}{2}(1-\frac{1}{p})-|\beta|/2}\} E_{|\beta|+[m/2]+1},$$

with $p \in [2, +\infty]$.

Let us state our global decay estimate for $u$.

**Theorem 5.4.** *Let $u(t)$ be a smooth global solution to problem (5.1), (5.2). Let $E_s = \max\{\|u(0)\|_{L^1}, \|u(0)\|_{H^s}\}$, and assume $E_{[m/2]+2}$ small enough. Let $p \in [2, +\infty]$. The following decay estimate holds*

$$\text{(5.44)} \quad \|D^\beta u(t)\|_{L^p} \leq C \min\{1, t^{-\frac{m}{2}(1-\frac{1}{p})-|\beta|/2}\} E_{|\beta|+[m/2]+1},$$

*with $C = C(E_{|\beta|+\sigma})$, for $\sigma$ large enough.*

*Remark* 5.5. For $m = 1$ we can estimate also the $L^1$ norm, since with the same computation as above we have

$$\|D^\beta u(t)\|_{L^1} \leq C\min\{1, t^{-\beta/2}\}\|u(0)\|_{L^1} + Ce^{-ct}\|D^\beta u(0)\|_{L^1}$$

$$+CE_1^2 \min\{1, t^{-\beta/2-1/2}\} \int_0^{t/2} \min\{1, s^{-1/2}\}ds$$

$$+CE_\beta^2 \min\{1, t^{-1/2-\beta/2}\} \int_{t/2}^t \min\{1, (t-s)^{-1/2}\}ds$$

$$+CE_{\beta+1}E_\beta \int_0^t e^{-c(t-s)} \min\{1, s^{-1/2-\beta/2}\}ds,$$

which yields

(5.45) $$\|D^\beta u(t)\|_{L^1} \leq C\left(\min\{1, t^{-\beta/2}\}E_{\beta+1} + e^{-ct}\|D^\beta u(0)\|_{L^1}\right),$$

with the constant $C = C(E_{\beta+2})$, so that, for $m = 1$, Theorem 5.4 holds for $p \in [1, +\infty]$.

5.2. **Decay estimates for the dissipative variables.** We now study the faster decay of the dissipative variables. Set $u_c = L_0 u(t)$ for the conservative variables, and $u_d(t) = L_- u(t)$ for the dissipative ones, where the projectors $L_0$ and $L_-$ are given by (4.14). We have that

$$D^\beta u_d(t) = D^\beta L_- K(t)u(0) + L_- \mathcal{K}(t)D^\beta u(0)$$

$$+ \sum_{\alpha=1}^m \int_0^{t/2} D^\beta D_{x_\alpha} L_- K(t-s)\bigl(Df_\alpha(0)u(s) - f_\alpha(u(s))\bigr)ds$$

$$+ \sum_{\alpha=1}^m \int_{t/2}^t D_{x_\alpha} L_- K(t-s)D^\beta\bigl(Df_\alpha(0)u(s) - f_\alpha(u(s))\bigr)ds$$

$$+ \int_0^{t/2} D^\beta L_- K(t-s)R_- L_-\bigl(g(u(s)) - Dg(0)u(s)\bigr)ds$$

$$+ \int_{t/2}^t L_- K(t-s)R_- L_- D^\beta\bigl(g(u(s)) - Dg(0)u(s)\bigr)ds$$

(5.46) $$+ \sum_{\alpha=1}^m \int_0^t L_- \mathcal{K}(t)D^\beta\bigl(D_{x_\alpha}\bigl(Df_\alpha(0)u(s) - f_\alpha(u(s))\bigr) + \bigl(g(u(s)) - Dg(0)u(s)\bigr)\bigr)ds.$$

As we see form the above formula, in this case one gains $t^{-1/2}$ in the estimates of the convolution with the smoothing kernels, because the principal terms in the initial data is $K_{-0}$ and in the convolutions are $DK_{-0}$ and $K_{--}(t)$, respectively, but no gain in the singular part $L_- \mathcal{K}(t)$.

We start with the $L^2$ norm of the $\beta$ derivative: we have

$$\|D^\beta u_d(t)\|_{L^2} \leq C\min\{1, t^{-m/4-|\beta|/2-1/2}\}\|u(0)\|_{L^1} + Ce^{-ct}\|u(0)\|_{H^{|\beta|}}$$

$$+C\int_0^{t/2} \min\{1, (t-s)^{-m/4-1-|\beta|/2}\}\|u(s)^2\|_{L^1}ds$$

(5.47) $$+C\int_{t/2}^t \min\{1, (t-s)^{-m/4-1}\}\Bigl(\sum_{\alpha=1}^m \|D^\beta u^2 h_\alpha(u(s))\|_{L^1} + \|D^\beta u^2 h(u(s))\|_{L^1}\Bigr)ds$$

$$+C\int_0^t e^{-c(t-s)}\Bigl(\sum_{\alpha=1}^m \|D^\beta D_{x_\alpha} u^2 h_\alpha(u(s))\|_{L^2} + \|D^\beta u^2 h(u(s))\|_{L^2}\Bigr)ds.$$

Using (5.36) and (5.33), we have

$$\sum_{\alpha=1}^{m} \|D^\beta u^2 h_\alpha(u(s))\|_{L^1} + \|D^\beta u^2 h(u(s))\|_{L^1} \leq C \min\{1, t^{-m/2-|\beta|/2}\}\left(E^2_{|\beta|+[m/2]} + E_{[m/2]+1}E_{|\beta|+[m/2]+1}\right). \quad (5.48)$$

Next, using (5.23) and then (5.44) and (5.40), yields

$$\sum_{\alpha=1}^{m} \|D^\beta D_{x_\alpha} u^2 h_\alpha(u(s))\|_{L^2} + \|D^\beta u^2 h(u(s))\|_{L^2} \leq C\|u\|_{L^\infty}\|D^\beta u\|_{H^{[m/2]+2}} \quad (5.49)$$

$$\leq C \min\{1, t^{-m/2}\}E_{[m/2]+1} \min\{1, t^{-m/4-|\beta|/2}\}E_{|\beta|+[m/2]+2}.$$

Therefore, we use the above inequalities in (5.47), to give

$$\|D^\beta u_d(t)\|_{L^2} \leq C \min\{1, t^{-m/4-|\beta|/2-1/2}\}E_{|\beta|}$$

$$+ C \min\{1, t^{-m/4-1-|\beta|/2}\}\left(\int_0^{t/2} \min\{1, s^{-m/2}\}ds\right) E^2_{[m/2]+1}$$

$$(5.50)$$

$$+ C \min\{1, t^{-m/2-|\beta|/2}\}\left(\int_{t/2}^{t} \min\{1, (t-s)^{-m/4-1}\}ds\right)\left(E^2_{|\beta|+[m/2]} + E_{[m/2]+1}E_{|\beta|+[m/2]+1}\right)$$

$$+ C\left(\int_0^t e^{-c(t-s)} \min\{1, s^{-3m/4-|\beta|/2}\}ds\right) E_{[m/2]+1}E_{|\beta|+[m/2]+2}.$$

Then, for $m \geq 2$, we obtain

$$\|D^\beta u_d(t)\|_{L^2} \leq C \min\{1, t^{-m/4-|\beta|/2-1/2}\}(E_{|\beta|} + E^2_{[m/2]+1} + E^2_{|\beta|+[m/2]} + E_{[m/2]+1}E_{|\beta|+[m/2]+2}),$$

which implies

$$\|D^\beta u_d(t)\|_{L^2} \leq C \min\{1, t^{-m/4-|\beta|/2-1/2}\}E_{[m/2]+|\beta|+1}, \quad (5.51)$$

with $C = C(E_{|\beta|+[m/2]+2})$.

About the $L^\infty$ norm, we can use (5.51) and, by arguing as in (5.41), we have

$$\|D^\beta u_d(t)\|_{L^\infty} \leq C \min\{1, t^{-m/2-|\beta|/2-1/2}\}E_{|\beta|+m+1}, \quad (5.52)$$

for $m \geq 2$. Actually, it is also possible as before, to show by a direct calculation, that the following estimate holds:

$$\|D^\beta u_d(t)\|_{L^\infty} \leq C \min\{1, t^{-m/2-|\beta|/2-1/2}\}E_{|\beta|+[m/2]+1}. \quad (5.53)$$

We also obtain the $L^p$-decay estimate

$$\|D^\beta u_d(t)\|_{L^p} \leq C \min\{1, t^{-\frac{m}{2}(1-\frac{1}{p})-|\beta|/2-1/2}\}E_{|\beta|+[m/2]+1}, \quad (5.54)$$

with $p \in [2+\infty]$, $m \geq 2$.

For $m = 1$, we recall that, thanks to (3.57), $L_-K(t)R_- = K_{--} = D^2 S_{--}$, where $S_{--}$ is a heat like kernel. Therefore, we will consider the equation

$$D^\beta u_d(t) = \left(D^\beta L_- K(t) + D^\beta L_- R(t)\right) u(0) + \mathcal{K}(t) D^\beta u(0)$$

$$+ \int_0^{t/2} D^{\beta+1}\left(L_- K(t-s) + L_- R(t-s)\right)\left(Df(0)u(s) - f(u(s))\right) ds$$

$$+ \int_{t/2}^t \left(L_- K(t-s) + L_- R(t-s)\right) D^{\beta+1}\left(Df(0)u(s) - f(u(s))\right) ds$$

$$+ \int_0^{t/2} D^\beta \left(D^2 S_{--}(t-s) + R_{--}(t-s)\right) L_-\left(g(u(s)) - Dg(0)u(s)\right) ds$$

$$+ \int_{t/2}^t DS_{--}(t-s) L_- D^{\beta+1}\left(g(u(s)) - Dg(0)u(s)\right) ds$$

$$+ \int_{t/2}^t R_{--}(t-s) L_- D^\beta \left(g(u(s)) - Dg(0)u(s)\right) ds$$

(5.55)
$$+ \int_0^t L_- \mathcal{K}(t-s) D^\beta \left(D\left(Df(0)u(s) - f(u(s))\right) + \left(g(u(s)) - Dg(0)u(s)\right)\right) ds.$$

Thanks to this decomposition we can prove the $L^2$ estimate for $m = 1$. Using (5.48) and, for the fifth integral, the $L^1$-estimate of $R_{--}$ coupled with (5.23), we have

$$\|D^\beta u_d(t)\|_{L^2} \leq C \min\{1, t^{-1/4 - \beta/2 - 1/2}\} E_\beta$$

$$+ C \min\{1, t^{-1/4 - 1 - \beta/2}\} \left(\int_0^{t/2} \min\{1, s^{-1/2}\} ds\right) E_1^2$$

(5.56)
$$+ C \min\{1, t^{-1-\beta/2}\} \left(\int_{t/2}^t \min\{1, (t-s)^{-3/4}\} ds\right) E_{\beta+2} E_{\beta+1}$$

$$+ C \min\{1, t^{-3/4 - \beta/2}\} \left(\int_{t/2}^t \min\{1, (t-s)^{-3/2}\} ds\right) E_1 E_{\beta+1}$$

$$+ C \left(\int_0^t e^{-c(t-s)} \min\{1, s^{-3/4 - \beta/2}\} ds\right) E_1 E_{\beta+2}.$$

Therefore, we have (5.51) for $m = 1$. On the other hand, we can estimate the $L^\infty$-norm as follows. First we have

$$\|D^\beta u_d(t)\|_{L^\infty} \leq C \min\{1, t^{-1-\beta/2}\} \|u(0)\|_{L^1} + C e^{-ct} \|D^\beta u(0)\|_{L^\infty}$$

$$+ C \int_0^{t/2} \min\{1, (t-s)^{-3/2 - \beta/2}\} \|u(s)^2\|_{L^1} ds$$

$$+ C \int_{t/2}^t \min\{1, (t-s)^{-1/2}\} \|D^{\beta+1} u^2 h(u(s))\|_{L^\infty} ds$$

$$+ C \int_{t/2}^t \min\{1, (t-s)^{-3/2}\} \|D^\beta u^2 h(u(s))\|_{L^\infty} ds$$

$$+ C \int_0^t e^{-c(t-s)} \left(\|D^{\beta+1} u^2 h(u(s))\|_{L^\infty} + \|D^\beta u^2 h(u(s))\|_{L^\infty}\right) ds.$$

By using (5.42), we have

(5.57) $$\|D^\beta u^2 h(u(t))\|_{L^\infty} \leq C \min\{1, t^{-1-\beta/2}\} E_{\beta+1} E_{\sup(\beta, 1)},$$

which yields

$$\begin{aligned}\|D^\beta u_d(t)\|_{L^\infty} \leq\ & C\min\{1, t^{-1-\beta/2}\}\|u(0)\|_{L^1} + Ce^{-ct}E_{\beta+1} \\
& + CE_1^2 \min\{1, t^{-3/2-\beta/2}\} \int_0^{t/2} \min\{1, s^{-1/2}\}ds \\
& + CE_{\beta+2}E_{\beta+1} \min\{1, t^{-3/2-\beta/2}\} \int_{t/2}^t \min\{1, (t-s)^{-1/2}\}ds \\
& + CE_{\beta+1}^2 \min\{1, t^{-1-\beta/2}\} \int_{t/2}^t \min\{1, (t-s)^{-3/2}\}ds \\
& + CE_{\beta+2}E_{\beta+1}\min\{1, t^{-1-\beta/2}\}.\end{aligned}$$

Therefore, we conclude

$$\|D^\beta u_d(t)\|_{L^\infty} \leq C\min\{1, t^{-1-\beta/2}\}E_{\beta+1}, \tag{5.58}$$

with $C = C(E_{\beta+2})$.

**Theorem 5.6.** *Under the assumptions of Theorem 5.4, we have the following decay estimates for the dissipative part of $u$:*

$$\|D^\beta u_d(t)\|_{L^p} \leq C\min\{1, t^{-\frac{m}{2}(1-\frac{1}{p})-1/2-|\beta|/2}\}E_{|\beta|+[m/2]+1}, \tag{5.59}$$

*with $C = C(E_{|\beta|+\sigma})$, for $\sigma$ large enough, and $p \in [2, +\infty]$.*

*Remark* 5.7. As previously, for $m = 1$, we can estimate the $L^1$ norm. By the decomposition (5.55), and using (5.48), we have

$$\begin{aligned}\|D^\beta u_d(t)\|_{L^1} \leq\ & C\min\{1, t^{-1/2-\beta/2}\}\|u(0)\|_{L^1} + Ce^{-ct}\|D^\beta u(0)\|_{L^1} \\
& + CE_1^2 \min\{1, t^{-1-\beta/2}\} \int_0^{t/2} \min\{1, s^{-1/2}\}ds \\
& + CE_{\beta+1}E_{\beta+2}\min\{1, t^{-1-\beta/2}\} \int_{t/2}^t \min\{1, (t-s)^{-1/2}\}ds \\
& + CE_{\beta+1}^2 \min\{1, t^{-1/2-\beta/2}\} \int_{t/2}^t \min\{1, (t-s)^{-3/2}\}ds \\
& + CE_{\beta+1}E_{\beta+2} \int_0^t e^{-c(t-s)}\min\{1, s^{-1/2-\beta/2}ds\}\end{aligned}$$

$$\leq C\min\{1, t^{-1/2-\beta/2}\}E_{\bar{b}+1}, \tag{5.60}$$

which yields

$$\|D^\beta u_d(t)\|_{L^1} \leq C\left(\min\{1, t^{-1/2-\beta/2}\}E_{\beta+1} + e^{-ct}\|D^\beta u(0)\|_{L^1}\right), \tag{5.61}$$

with the constant $C = C(E_{\beta+2})$, so that, for $m = 1$, Theorem 5.6 holds for $p \in [1, +\infty]$.

5.3. **Decay estimates for the time derivative.** We estimate now the decay of the time derivative of the solution. Directly from equations we obtain

$$\|D^\beta u_t(t)\|_{L^p} \leq \sum_{\alpha=1}^m \|D^\beta D_{x_\alpha} f_{c,\alpha}(u)\|_{L^p} + C\|D^\beta u_d\|_{L^p} + \|D^\beta(g(u) - Dg(0)u)\|_{L^p}$$

$$\leq C\min\{1, t^{-\frac{m}{2}(1-\frac{1}{p})-|\beta|/2-1/2}\}E_{|\beta|+[m/2]+2}, \tag{5.62}$$

where $p \in [2, +\infty]$. For $m = 1$, as previously,

$$\|D^\beta u_t(t)\|_{L^1} \leq C\min\{1, t^{-1/2-\beta/2}\}E_{\beta+2} + Ce^{-ct}\|D^{\beta+1}u(0)\|_{L^1}. \tag{5.63}$$

About the dissipative variables, we write

$$(u_t)_t + \sum_{\alpha=1}^m Df_\alpha(0)(u_t)_{x_\alpha} - Dg(0)u_t = \sum_{\alpha=1}^m ((Df_\alpha(0) - Df_\alpha(u))u_t)_{x_\alpha} + (Dg(u) - Dg(0))u_t,$$

so that we obtain

$$D^\beta u_{d,t}(t) = D^\beta L_- K(t)u_t(0) + \mathcal{K}(t)D^\beta u_t(0)$$

$$+ \sum_{\alpha=1}^m \int_0^{t/2} D^\beta D_{x_\alpha} L_- K(t-s)\Big((Df_\alpha(0)u(s) - Df_\alpha(u))u_t(s)\Big)ds$$

$$+ \sum_{\alpha=1}^m \int_{t/2}^t D_{x_\alpha} L_- K(t-s)D^\beta\Big((Df_\alpha(0) - Df_\alpha(u))u_t(s)\Big)ds$$

$$+ \int_0^{t/2} D^\beta L_- K(t-s)R_- L_-\Big((Dg(u) - Dg(0))u_t(s)\Big)ds$$

$$+ \int_{t/2}^t L_- K(t-s)R_- D^\beta L_-\Big((Dg(u) - Dg(0))u_t(s)\Big)ds$$

(5.64)
$$+ \sum_{\alpha=1}^m \int_0^t \mathcal{K}(t)D^\beta\Big(D_{x_\alpha}((Df_\alpha(0) - Df_\alpha(u))u_t(s)) + (Dg(u) - Dg(0))u_t(s)\Big)ds.$$

It holds

$$L_- D^\beta K(t)R_0 L_0 u_t(0) = -\sum_\alpha D^\beta D_{x_\alpha} L_- K(t)R_0 L_0 f_\alpha(u(0)),$$

and

$$L_- D^\beta K(t)R_- L_- u_t(0) = -\sum_\alpha D^\beta D_{x_\alpha} L_- K(t)R_- L_- f_\alpha(u(0)) + D^\beta L_- K(t)R_- L_- g(u(0)).$$

Since $f_\alpha(0) = g(0) = 0$, we have $\|f_\alpha(u)\|_{L^1}, \|g(u)\|_{L^1} = O(1)\|u\|_{L^1}$, which can be used, for $m \geq 2$, to yield

$$\|D^\beta u_{d,t}(t)\|_{L^2} \leq C\min\{1, t^{-m/4-|\beta|/2-1}\}\|u(0)\|_{L^1} + Ce^{-ct}E_{|\beta|+[m/2]+2}$$

$$+ CE_{[m/2]+2}\min\{1, t^{-m/4-1-|\beta|/2}\}\int_0^{t/2} \min\{1, s^{-m/2-1/2}\}ds$$

$$+ CE_{|\beta|+[m/2]+2}\min\{1, t^{-m/2-|\beta|/2-1/2}\}\int_{t/2}^t \min\{1, (t-s)^{-m/4-1}\}ds$$

$$+ CE_{|\beta|+[m/2]+3}\min\{1, t^{-3m/4-|\beta|/2-1/2}\}$$

(5.65)
$$\leq C\min\{1, t^{-m/4-|\beta|/2-1}\}E_{|\beta|+[m/2]+3}.$$

It follows by using the same analysis of (5.58) that for all $p \in [2, +\infty]$

(5.66)
$$\|D^\beta u_{d,t}(t)\|_{L^p} \leq C\min\{1, t^{-\frac{m}{2}(1-\frac{1}{p})-|\beta|/2-1}\}E_{|\beta|+[m/2]+3}.$$

Repeating the computations we did above for $m = 1$, it follows that the above estimate holds also for $m = 1$.

**Theorem 5.8.** *Under the assumptions of Theorem 5.4, we have the following decay estimates for $u_t$*

(5.67)
$$\|D^\beta u_t(t)\|_{L^p} \leq C\min\{1, t^{-\frac{m}{2}(1-\frac{1}{p})-1/2-|\beta|/2}\}E_{|\beta|+[m/2]+2},$$

(5.68)
$$\|D^\beta u_{d,t}(t)\|_{L^p} \leq C\min\{1, t^{-\frac{m}{2}(1-\frac{1}{p})-1-|\beta|/2}\}E_{|\beta|+[m/2]+3},$$

*where $C = C(E_{|\beta|+\sigma})$ for $\sigma$ large enough and $p \in [2, \infty]$.*

For $m = 1$, as previously,

(5.69)
$$\|D^\beta u_{d,t}(t)\|_{L^1} \leq C\min\{1, t^{-1-\beta/2}\}E_{\beta+3} + Ce^{-ct}\|D^{\beta+1}u(0)\|_{L^1}.$$

5.4. **Decay to linear solution.** We consider here the difference among the solution of the nonlinear equation (5.1) and the linearized one

$$u_t + \sum_{\alpha=1}^{m} Df_\alpha(0)u_{x_\alpha} = \begin{pmatrix} 0 \\ D_{u_d}q(0)u_d \end{pmatrix}, \tag{5.70}$$

where we have already considered the conservative-dissipative variable pair. The idea is that if the dimension $m$ is sufficiently large, then the decay of the non linear parts is faster than the linear part. Since it is easy to show that the following results do not hold in the case $m = 1$, we will consider in the following $m \geq 2$, thus estimating only the $L^p$ norm, $p \in [2, +\infty]$.

By using the representation (5.7), it follows that

$$u(t) - u_l(t) = \sum_{\alpha=1}^{m} \int_0^t D_{x_\alpha}\Gamma(t-s)\big(Df_\alpha(0)u(s) - f_\alpha(u(s))\big)ds$$

$$+ \int_0^t \Gamma(t-s)(g(u(s)) - Dg(0)u(s))ds. \tag{5.71}$$

Repeating the estimates leading to (5.38) and (5.51), it follows that for $m \geq 3$

$$\|D^\beta(u(t) - u_l(t))\|_{L^2} \leq \min\big\{1, t^{-\frac{m}{4}-|\beta|/2-1/2}\big\}E^2_{|\beta|+[m/2]+1}. \tag{5.72}$$

By arguing as previously, it follows also that

$$\|D^\beta(u(t) - u_l(t))\|_{L^p} \leq \min\big\{1, t^{-\frac{m}{2}(1-\frac{1}{p})-|\beta|/2-1/2}\big\}E^2_{|\beta|+[m/2]+1}, \tag{5.73}$$

for $m \geq 3$. If one tries to repeat the above computations for $m = 2$, one finds that there is a critical integral of the form

$$\mathcal{I} = \int_0^t L_-K(t-s)R_-L_-(g(u(s)) - Dg(0)u(s))ds,$$

which we can only estimate at order $\ln t/t$, since, directly using Theorem 4.2, we obtain that

$$\|\mathcal{I}\|_{L^2} \leq C \int_0^t \min\big\{1, (t-s)^{-1}\big\}\min\big\{1, s^{-1}\big\}ds.$$

However, using Remark 4.6, we can actually write

$$\mathcal{I} = O(1)\sum_\alpha \int_0^t \left(\int_0^{t-s} L_-K(t-s-\tau)A_\alpha R_-e^{\tau D_{u_d}q(0)}d\tau\right)D_{x_\alpha}(g(u(s)) - Dg(0)u(s))ds.$$

Therefore, estimating the Kernel in $L^1$ and the term $D_{x_\alpha}(g(u(s)) - Dg(0)u(s))$ in $L^2$, we obtain

$$\|\mathcal{I}\|_{L^2} \leq C \int_0^t \min\big\{1, (t-s)^{-1/2}\big\}\min\big\{1, s^{-3/2}\big\} \leq C\min\big\{1, s^{-1}\big\}.$$

We thus conclude with the following result.

**Theorem 5.9.** *Let $u_l$ be the solution of problem (5.70), (5.2), under the assumptions of Theorem 5.4, for $m \geq 2$ and $p \in [2, \infty]$, we have the following decay estimate*

$$\|D^\beta(u(t) - u_l(t))\|_{L^p} \leq C\min\big\{1, t^{-\frac{m}{2}(1-\frac{1}{p})-|\beta|/2-1/2}\big\}E_{|\beta|+[m/2]+1}, \tag{5.74}$$

*with $C = C(E_{|\beta|+\sigma})$, for $\sigma$ large enough.*

We notice now that for the conservative part we have that

$$L_0u_l(t) = L_0\Gamma(t)(R_0L_0u(0) + R_-L_-u(0))$$

$$= K_{00}(t)L_0u(0) + L_0K(t)R_-L_-u(0) + L_0\mathcal{K}(t)u(0).$$

So, using Theorem 4.2, and by arguing as for Theorem 5.9, we obtain another interesting approximation.

**Theorem 5.10.** *Under the assumptions of Theorem 5.9, the following decay estimate holds*

$$\|D^\beta(u_c(t) - K_{00}(t)L_0u(0))\|_{L^p} \leq C\min\big\{1, t^{-\frac{m}{2}(1-\frac{1}{p})-|\beta|/2-1/2}\big\}E_{|\beta|+[m/2]+1}. \tag{5.75}$$

Let us notice that, by definition, the pseudo-differential operator $K_{00}$ is always fully parabolic.

### 5.5. Chapman-Enskog expansion.

We show now how the solutions to the parabolic Chapman-Enskog expansion approximate the conservative part of the solutions to the nonlinear hyperbolic problem:

$$(5.76) \qquad u_t + \sum_{\alpha=1}^m A_\alpha(u) u_{x_\alpha} = \begin{pmatrix} 0 \\ q(u) \end{pmatrix},$$

where $A_\alpha(u) = Df_\alpha(u)$. We use the conservative-dissipative decomposition of $u$:

$$(5.77) \qquad u_{c,t} + \sum_{\alpha=1}^m A_{\alpha,11}(0) u_{c,x_\alpha} + \sum_{\alpha=1}^m A_{\alpha,12}(0) u_{d,x_\alpha} = L_0 \sum_{\alpha=1}^m (A_\alpha(0) u - f_\alpha(u))_{x_\alpha} ;$$

$$(5.78) \quad u_{d,t} + \sum_{\alpha=1}^m A_{\alpha,21}(0) u_{c,x_\alpha} + \sum_{\alpha=1}^m A_{\alpha,22}(0) u_{d,x_\alpha} = D_{u_d} q(0) u_d + L_- \sum_{\alpha=1}^m (A_\alpha(0) u - f_\alpha(u))_{x_\alpha} + (q(u) - D_{u_d} q(0) u_d).$$

We can compute $u_d$ using (5.78), which yields, inserting in (5.77):

$$(5.79) \qquad \begin{aligned} & u_{c,t} + \sum_{\alpha=1}^m A_{\alpha,11}(0) u_{c,x_\alpha} + \sum_{\beta=1}^m \sum_{\alpha=1}^m A_{\alpha,12}(0) (D_{u_d} q(0))^{-1} A_{\beta,21}(0) u_{c,x_\alpha x_\beta} \\ &= L_0 \sum_{\alpha=1}^m (A_\alpha(0) u - f_\alpha(u))_{x_\alpha} + \sum_{\alpha=1}^m A_{\alpha,12}(0) (D_{u_d} q(0))^{-1} (q(u) - D_{u_d} q(0) u_d)_{x_\alpha} \\ &\quad - \sum_{\alpha=1}^m A_{\alpha,12}(0) (D_{u_d} q(0))^{-1} \left( u_{d,t x_\alpha} + \sum_{\beta=1}^m A_{\beta,22}(0) u_{d,x_\alpha x_\beta} - \sum_{\beta=1}^m L_- \left( A_\beta(0) u - f_\beta(u) \right)_{x_\alpha x_\beta} \right). \end{aligned}$$

We consider the linear parabolic equation

$$(5.80) \qquad w_t + \sum_{\alpha=1}^m A_{\alpha,11}(0) w_{x_\alpha} + \sum_{\beta=1}^m \sum_{\alpha=1}^m A_{\alpha,12}(0) (D_{u_d} q(0))^{-1} A_{\beta,21}(0) w_{x_\alpha x_\beta} = 0,$$

and we denote by $u_p(t)$ the solution of the weakly parabolic equation (5.80) with

$$(5.81) \qquad u_p(0) = L_0 u(0).$$

Using Remark 4.3 and by arguing again as for Theorem 5.9, it is possible to prove the following result.

**Theorem 5.11.** *Let $u_p$ be the solution of problem* (5.80), (5.81), *under the assumptions of Theorem 5.4, for $m \geq 2$ and $p \in [2, \infty]$, we have the following decay estimate*

$$(5.82) \qquad \|D^\beta(u_c(t) - u_p(t))\|_{L^p} \leq C \min\left\{1, t^{-\frac{m}{2}(1-\frac{1}{p}) - |\beta|/2 - 1/2}\right\} E_{|\beta| + [m/2] + 1},$$

*with $C = C(E_{|\beta| + \sigma})$, for $\sigma$ large enough.*

Notice that the same faster decay holds for the difference between the solution $u_p$ to the weakly parabolic problem and the solution $K_{00} L_0 u(0)$ of the "fully" parabolic pseudo-differential problem.

*Example* 5.12. **Rotationally invariant systems** Consider the isentropic dissipative Euler equations

$$(5.83) \qquad \begin{cases} \rho_t + \mathrm{div}(\rho v) = 0, \\ (\rho v)_t + \mathrm{div}(\rho v \otimes v) + \frac{1}{\gamma} \nabla \rho^\gamma = -v. \end{cases}$$

We can linearize the system around the constant state $(\bar\rho, \bar v) = (1, 0)$, so obtaining system (4.40) of Example 4.7. In that case we can immediately apply Theorems 5.4, 5.6, 5.9, and 5.11. In particular, by eliminating $v$ in (4.40), we obtain the estimate

$$\|D^\beta(\rho(t) - \rho_w(t))\|_{L^p} + \|D^\beta(\rho(t) - \rho_p(t))\|_{L^p} \leq C \min\left\{1, t^{-\frac{m}{2}(1-\frac{1}{p}) - |\beta|/2 - 1/2}\right\},$$

where $\rho_w$ and $\rho_p$ are respectively the solutions of the $m$-dimensional dissipative wave equation equation

$$\rho_{w,t} + \rho_{w,tt} - \Delta \rho_w = 0,$$

and the $m$-dimensional heat equation

$$\rho_{p,t} - \Delta \rho_p = 0.$$

These estimates improve on previous results about this problem contained in [34] and [7].

Consider now the relaxation system

(5.84)
$$\begin{cases} \rho_t + \text{div}(\rho v) = 0, \\ (\rho v)_t + \text{div}(\rho R) + \nabla \rho = 0, \\ (\rho R)_t + \nabla(\rho v) = \rho v \otimes v - \rho R. \end{cases}$$

Its local relaxation limit is given by the (non dissipative) isentropic Euler equations. However, its linearized version around the state $(\bar\rho, \bar v, \bar R) = (1, 0, 0)$, is just given by system (4.43) of Example 4.7. Again, we can explicitly identify the asymptotic limits, with analogous decay rates, in terms of the linear hyperbolic system (4.43) and, thanks to the analysis in Example 4.7, of the fully parabolic system (4.49), which corresponds to the kernel $K_{00}$ given by (4.44). Finally, thanks to Theorem 5.11, the same behavior is shown by its Chapman-Enskog expansion, which is given in this case by the weakly parabolic system

$$\begin{cases} \rho_t + \text{div} v = 0, \\ v_t + \nabla \rho = \Delta v. \end{cases}$$

*The Chapman-Enskog expansion in the case $m = 1$.* For $m = 1$, we need to consider together with the linear part the nonlinear terms of the order of $u^2$, because the decay of $u^2$ convoluted with the linear kernel and integrated in time gives the same decay estimate of $u$. We will prove that for all $0 \leq \mu < 1/2$, the difference among the conservative variables and the solution to an approximated Chapman-Enskog expansion decays as $t^{-1/2(1-1/p)-\mu}$ in $L^p$ if the initial data is sufficiently small: their size goes to 0 as $\mu \to 1/2$.

Let us introduce the operators

(5.85)
$$\tilde A = \frac{1}{2}\left(L_0 D^2_{u_c} f(0) - A_{12}(0)(D_{u_d} q(0))^{-1} D^2_{u_c} q(0)\right),$$

(5.86)
$$\tilde B = A_{12}(0)(D_{u_d} q(0))^{-1} A_{21}(0).$$

We rewrite (5.79) as

(5.87)
$$u_{c,t} + \left(A_{11}(0) u_c + \tilde A(u_c, u_c)\right)_x + \tilde B u_{c,xx} = S_x,$$

with

$$\begin{aligned} S = & \ L_0\left(A(0)u - f(u) + \tfrac{1}{2} D^2_{u_c} f(0)(u_c, u_c)\right) \\ & + A_{12}(0)(D_{u_d} q(0))^{-1}\left(q(u) - D_{u_d} q(0) u_d - \tfrac{1}{2} D^2_{u_c} q(0)(u_c, u_c)\right) \\ & - A_{12}(0)(D_{u_d} q(0))^{-1}\left(u_{d,t} + A_{22}(0) u_{d,x} - L_-\left(A(0)u - f(u)\right)_x\right). \end{aligned}$$

In the same way, we replace (5.80) by the nonlinear parabolic equation

(5.88)
$$w_t + \left(A_{11}(0)w + \tilde A(w, w)\right)_x + \tilde B w_{xx} = 0.$$

We introduce $F_\beta$, with $F_1 = E_1$ and, if $\beta \geq 1$,

(5.89)
$$F_{\beta+1} = \begin{cases} E_{\beta+1}, & \text{if } p \in [2, \infty], \\ E_{\beta+1} + \|D^\beta u(0)\|_{L^1}, & \text{otherwise.} \end{cases}$$

**Theorem 5.13.** *Let $u_p$ be the solution of problem (5.88), (5.81), under the assumptions of Theorem 5.4, for $m = 1$ and $p \in [1, \infty]$, for $\mu \in [0, 1/2)$, if $E_1$ sufficiently small with respect to $(1/2 - \mu)$, then we have the following decay estimate*

(5.90)
$$\|D^\beta(u_c(t) - u_p(t))\|_{L^p} \leq C \min\left\{1, t^{-\frac{1}{2}(1-\frac{1}{p})-\mu-\beta/2}\right\} F_{\beta+4},$$

*where $C = C(\mu, F_{\beta+\sigma})$, for $\sigma$ large enough.*

*Proof.* We denote by $\Gamma_p(t)$ the Green kernel of the linear parabolic equation

$$w_t + A_{11}(0)w_x + \tilde{B}w_{xx} = 0. \tag{5.91}$$

Using Remark 4.3, $\Gamma_p(t)$ can be written as

$$\Gamma_p(t) = K_{00}(t) + \tilde{\mathcal{K}}(t) + \tilde{R}(t), \tag{5.92}$$

where $K_{00}(t)$ is the 00 component of the principal part $K(t)$ of the relaxation kernel $\Gamma(t)$, given by (3.57).

We take the difference among $u_c(t)$ and $u_p(t)$:

$$\begin{aligned}
D^\beta(u_c(t) - u_p(t)) &= \int_0^{t/2} D^\beta D\big(K_{00}(t-s) + \tilde{R}(t-s)\big)\big(\tilde{A}(u_p(s), u_p(s)) - \tilde{A}(u_c(s), u_c(s))\big)ds \\
&\quad + \int_0^{t/2} D^\beta D\big(K_{00}(t-s) + \tilde{R}(t-s)\big)S(s)ds \\
&\quad + \int_{t/2}^t D(K_{00}(t-s) + \tilde{R}(t-s))D^\beta\big(\tilde{A}(u_p(s), u_p(s)) - \tilde{A}(u_c(s), u_c(s)) + S(s)\big)ds \\
&\quad + \int_0^t \tilde{\mathcal{K}}(t-s)D^\beta D\big(\tilde{A}(u_p(s), u_p(s)) - \tilde{A}(u_c(s), u_c(s)) + S(s)\big)ds. 
\end{aligned} \tag{5.93}$$

By the previous estimates on $D^\beta u$, $D^\beta u_d$, $D^\beta u_{d,t}$, we have that

$$\left\|D^\beta S\right\|_{L^p} \leq C \min\{1, t^{-1/2(1-1/p)-1-\beta/2}\}F_{\beta+3}. \tag{5.94}$$

Let us define, for a fixed $\mu \in [0, 1/2)$,

$$m_0(t) \doteq \sup_{0 \leq \tau \leq t}\{\max\{1, \tau^{1/4+\mu}\}\|u_c(\tau) - u_p(\tau)\|_{L^2}\}. \tag{5.95}$$

Taking the $L^2$ norm of (5.93), for $\beta = 0$ we have

$$\begin{aligned}
\|u_c(t) - u_p(t)\|_{L^2} &\leq CE_1 m_0(t) \int_0^t \min\{1, (t-s)^{-3/4}\}\min\{1, s^{-1/2-\mu}\}ds \\
&\quad + CF_3 \int_0^t \min\{1, (t-s)^{-3/4}\}\min\{1, s^{-1}\}ds \\
&\quad + C(E_2 E_1 + F_4)\int_0^t e^{-c(t-s)}\min\{1, s^{-5/4}\}ds \\
&\leq C\min\{1, s^{-1/4-\mu}\}(E_1 E_2 + E_1 + F_4 + (1/2 - \mu)^{-1}E_1 m_0(t)).
\end{aligned}$$

It follows thus

$$m_0(t) \leq CF_4, \tag{5.96}$$

for $E_1$ sufficiently small with respect to $(1/2 - \mu)$.

Assume now that for $\gamma < \beta$

$$\|D^\gamma(u_c - u_p)(t)\|_{L^2} \leq C(\mu)\min\{1, t^{-1/4-\mu-\gamma/2}\}F_{\gamma+4}, \tag{5.97}$$

with $\mu < 1/2$, and set

$$m_\beta(t) \doteq \sup_{0 \leq \tau \leq t}\{\max\{1, \tau^{1/4+\mu+\beta/2}\}\|D^\beta(u_c(\tau) - u_p(\tau))\|_{L^2}\}. \tag{5.98}$$

Using the induction assumption (5.97), we obtain

$$\begin{aligned}
\|D^\beta(\tilde{A}(u_c(s), u_c(s)) - \tilde{A}(u_p(s), u_p(s)))\|_{L^1} &\leq C\sum_{\alpha=0}^{\beta-1}\left(\|D^{\beta-\alpha}u_c\|_{L^2} + \|D^{\beta-\alpha}u_p\|_{L^2}\right)\|D^\alpha(u_c - u_p)\|_{L^2} \\
&\quad + C(\|u_c\|_{L^2} + \|u_p\|_{L^2})\|D^\beta(u_c - u_p)\|_{L^2} \\
&\leq C\min\{1, t^{-1/2-\mu-\beta/2}\}(C(\mu)E_{\beta+1}F_{\beta+3} + E_1 m_\beta(t)).
\end{aligned}$$

Using this inequality, (5.94) and (5.96) in (5.93) yields

$$\|D^\beta(u_c(t) - u_p(t))\|_{L^2} \le CE_1^2 \int_0^{t/2} \min\{1, (t-s)^{-3/4-\beta/2}\} \min\{1, s^{-1/2-\mu}\} ds$$

$$+ CF_3 \int_0^{t/2} \min\{1, (t-s)^{-3/4-\beta/2}\} \min\{1, s^{-1}\} ds$$

$$+ C(C(\mu)E_{\beta+1}F_{\beta+3} + E_1 m_\beta(t)) \int_{t/2}^t \min\{1, (t-s)^{-3/4}\} \min\{1, s^{-1/2-\mu-\beta/2}\} ds$$

$$+ CF_{\beta+4} \int_{t/2}^t \min\{1, (t-s)^{-3/4}\} \min\{1, s^{-1-\beta/2}\} ds$$

$$+ C(E_{\beta+2}E_1 + F_{\beta+4}) \int_0^t e^{-c(t-s)} \min\{1, s^{-5/4-\beta/2}\} ds.$$

It follows that $m_\beta(t) \le C(\mu) F_{\beta+4}$, and we have (5.97) for $\beta$.

As for the estimate (5.41), we have the $L^\infty$-estimate

(5.99) $$\|D^\beta(u_c(t) - u_p(t))\|_{L^\infty} \le C(\mu) \min\{1, t^{-1/2-\mu-\beta/2}\} F_{\beta+4}.$$

Finally, we estimate the $L^1$-norm in (5.93):

$$\|D^\beta(u_c(t) - u_p(t))\|_{L^1} \le CE_1^2 \int_0^{t/2} \min\{1, (t-s)^{-1/2-\beta/2}\} \min\{1, s^{-1/2-\mu}\} ds$$

$$+ CF_3 \int_0^{t/2} \min\{1, (t-s)^{-1/2-\beta/2}\} \min\{1, s^{-1}\} ds$$

$$+ C(C(\mu)E_{\beta+1}F_\beta + E_1 m_{\beta+3}(t)) \int_{t/2}^t \min\{1, (t-s)^{-1/2}\} \min\{1, s^{-1/2-\mu-\beta/2}\} ds$$

$$+ CF_{\beta+4} \int_{t/2}^t \min\{1, (t-s)^{-1/2}\} \min\{1, s^{-1-\beta/2}\} ds$$

$$+ \Big(C(\mu)(E_{\beta+2}F_{\beta+4} + E_1 F_{\beta+5}) + F_{\beta+4}\Big) \int_0^t e^{-c(t-s)} \min\{1, s^{-1-\beta/2}\} ds.$$

It follows that

(5.100) $$\|D^\beta(u_c(t) - u_p(t))\|_{L^1} \le C \min\{1, t^{-\mu-\beta/2}\} F_{\beta+4}.$$

Therefore, we obtain the conclusion. □

*Example* 5.14. **The $p$-system with relaxation.** We can apply Theorem 5.13 to the Example 2.10. In this case the Chapman-Enskog expansion is given by the semilinear parabolic equation

(5.101) $$u_{p,t} + h'(0)u_{p,x} + \frac{1}{2}h''(0)(u_p^2)_x - (\lambda^2 - a^2)u_{p,xx} = 0.$$

For previous results about this example see [6] and [22]. Notice that in [6], the data are chosen in a special class, which allows to take $\mu = 1/2$ in Theorem 5.13. However, even for this special example, our C-D decomposition gives a more precise description on the behavior of the solution, in terms of the dissipative part $u_d = (\lambda^2 - a^2)^{-\frac{1}{2}}(v - au)$.